\documentclass[11pt]{article}
\usepackage{latexsym,amssymb}
\usepackage{psfrag}
\usepackage{amsmath}
\usepackage{color}

\newtheorem{Theorem}{\bf Theorem}[section]
\newtheorem{Lemma}{\bf Lemma}[section]
\newtheorem{Proposition}{\bf Proposition}[section]
\newtheorem{Corollary}{\bf Corollary}[section]
\newtheorem{Remark}{\bf Remark}[section]
\newtheorem{Example}{\bf Example}[section]
\newtheorem{Definition}{\bf Definition}[section]

\newenvironment{theorem}{\begin{Theorem}$\!\!\!$}{\end{Theorem}}
\newenvironment{lemma}{\begin{Lemma}$\!\!\!$}{\end{Lemma}}

\newenvironment{corollary}{\begin{Corollary}$\!\!\!$}{\end{Corollary}}
\newenvironment{remark}{\begin{Remark}$\!\!\!$}{\end{Remark}}

\newenvironment{definition}{\begin{Definition}$\!\!\!$}{\end{Definition}}

\def\Xint#1{\mathchoice
{\XXint\displaystyle\textstyle{#1}}%
{\XXint\textstyle\scriptstyle{#1}}%
{\XXint\scriptstyle\scriptscriptstyle{#1}}%
{\XXint\scriptscriptstyle\scriptscriptstyle{#1}}%
\!\int}
\def\XXint#1#2#3{{\setbox0=\hbox{$#1{#2#3}{\int}$}
\vcenter{\hbox{$#2#3$}}\kern-.5\wd0}}

\def\dashint{\Xint-}

\numberwithin{equation}{section}

\begin{document}

\title{Solvability of the heat equation\\ with a nonlinear boundary condition}
\author{
Kotaro Hisa and Kazuhiro Ishige
\thanks{The second author of this paper was supported
by the Grant-in-Aid for Scientific Research (A)(No.~15H02058),
from Japan Society for the Promotion of Science.}
}
\date{}
\maketitle
\begin{abstract}
In this paper we obtain necessary conditions and sufficient conditions 
for the solvability of the problem
\begin{equation*}
({\rm P})\qquad
\left\{
\begin{array}{ll}
\partial_t u=\Delta u,\qquad & x\in{\bf R}^N_+,\,\,\,t>0,\vspace{3pt}\\
\partial_\nu u=u^p, & x\in\partial{\bf R}^N_+,\,\,\,t>0,\vspace{3pt}\\
u(x,0)=\mu(x)\ge 0,\qquad  & x\in D:=\overline{{\bf R}^N_+},
\end{array}
\right.
\qquad\qquad
\end{equation*}
where $N\ge 1$, $p>1$ and $\mu$ is a nonnegative measurable function in ${\bf R}^N_+$ 
or a Radon measure in ${\bf R}^N$ with $\mbox{supp}\,\mu\subset D$. 
Our sufficient conditions and necessary conditions enable us 
to identify the strongest singularity of the initial data for the solvability 
for problem~$({\rm P})$. 
Furthermore, as an application, 
we obtain optimal estimates of the life span of the minimal solution of $({\rm P})$
with $\mu=\kappa\varphi$ as $\kappa\to 0$ or $\kappa\to\infty$. 
\end{abstract}
\vspace{15pt}
\noindent Addresses:

\smallskip
\noindent K. H.:  Mathematical Institute, Tohoku University,
Aoba, Sendai 980-8578, Japan\\
\noindent 
E-mail: {\tt kotaro.hisa.s5@dc.tohoku.ac.jp}\\

\noindent K. I.: Mathematical Institute, Tohoku University,
Aoba, Sendai 980-8578, Japan\\
\noindent 
E-mail: {\tt ishige@m.tohoku.ac.jp}\\
\vspace{15pt}
\newpage
\section{Introduction}
We are interested in finding necessary conditions and sufficient conditions 
on the initial data for the solvability of problem
\begin{equation}
\label{eq:1.1}
\left\{
\begin{array}{ll}
\partial_t u=\Delta u,\qquad & x\in{\bf R}^N_+,\,\,\,t>0,\vspace{3pt}\\
\partial_\nu u=u^p, & x\in\partial{\bf R}^N_+,\,\,\,t>0,\vspace{3pt}
\end{array}
\right.
\end{equation}
with the initial condition
\begin{equation}
\label{eq:1.2}
u(x,0)=\mu(x)\ge 0,\qquad x\in D:=\overline{{\bf R}^N_+},
\end{equation}
where $\mu$ is a nonnegative measurable function in ${\bf R}^N_+$ 
or a Radon measure in ${\bf R}^N$ with $\mbox{supp}\,\mu\subset D$. 
For the solvability of problem~\eqref{eq:1.1} with \eqref{eq:1.2}, 
sufficient conditions have been studied in many papers 
(see e.g., \cite{ACR1}, \cite{AR}, \cite{DFL}, \cite{FiloK}, \cite{GL}, \cite{IK}, \cite{IS1} and \cite{IS2}). 
However little is known concerning necessary conditions
and the strongest singularity of initial data 
for which problem~\eqref{eq:1.1} possesses a local-in-time nonnegative solution is still open
as far as we know. 

In 1985, Baras and Pierre~\cite{BP} studied necessary conditions on the initial data
for the existence of nonnegative solutions of 
$$
\partial_t u=\Delta u+u^q,\quad x\in{\bf R}^N,\,\,t>0,
$$
where $N\ge 1$ and $q>1$. 
Recently, 
the authors of this paper \cite{HI} proved the existence and the uniqueness of the initial trace of a nonnegative solution of 
a fractional semilinear heat equation 
\begin{equation}
\label{eq:1.3}
\partial_t u=-(-\Delta)^{\frac{\theta}{2}}u+u^q,\quad  x\in{\bf R}^N,\,\,t>0,
\end{equation}
where $0<\theta\le 2$ and $q>1$. 
Furthermore, they showed that, 
if problem~\eqref{eq:1.3} possesses a local-in-time nonnegative solution, 
then its initial trace $\mu$ satisfies the following:
\begin{itemize}
  \item[{\rm (1)}] 
  $\displaystyle{\sup_{x\in{\bf R}^N}\mu(B(x,T^\frac{1}{\theta}))\le \gamma_1\,T^{\frac{N}{\theta}-\frac{1}{q-1}}}$ 
  if $1<q<q_\theta$;
  \item[{\rm (2)}] 
  $\displaystyle{\sup_{x\in{\bf R}^N}\mu(B(x,\sigma))\le \gamma_1\,
  \biggr[\log\biggr(e+\frac{T^{\frac{1}{\theta}}}{\sigma}\biggr)\biggr]^{-\frac{N}{\theta}}}$ 
  for all $0<\sigma<T^{\frac{1}{\theta}}$ if $q=q_\theta$; 
  \item[{\rm (3)}] 
  $\displaystyle{\sup_{x\in{\bf R}^N}\mu(B(x,\sigma))\le \gamma_1\,
  \sigma^{N-\frac{\theta}{p-1}}}$ 
  for all $0<\sigma<T^{\frac{1}{\theta}}$ if $q>q_\theta$.
\end{itemize} 
Here $q_\theta:=1+\theta/N$. 
In \cite{HI}, developing the arguments in \cite{IKS} and \cite{RS}, 
they also obtained sufficient conditions on the initial data for the existence of the solution of \eqref{eq:1.3} 
and identified the strongest singularity of the initial data for which 
the Cauchy problem to \eqref{eq:1.3} possesses a local-in-time nonnegative solution. 

In this paper, motivated by \cite{HI}, 
we show the existence and the uniqueness of the initial trace 
of a nonnegative solution of \eqref{eq:1.1} and obtain necessary conditions 
on the existence of nonnegative solutions of \eqref{eq:1.1} and \eqref{eq:1.2}. 
We also obtain new sufficient conditions on 
the existence of nonnegative solutions of \eqref{eq:1.1} and \eqref{eq:1.2}. 
Our necessary conditions and sufficient conditions enable us to 
identify the strongest singularity of initial data for which problem~\eqref{eq:1.1} 
possesses a local-in-time nonnegative solution. 
Surprisingly, the strongest singularity depends on whether it exists on $\partial{\bf R}^N_+$ or not 
(see Corollary~\ref{Corollary:1.1} and Section~6). 
Furthermore, 
we study how the life span of the solution of problem~\eqref{eq:1.1} with \eqref{eq:1.2} 
depends on the behavior of the initial data near the boundary and at the space infinity. 
See Section~6.
\vspace{5pt}

We introduce some notation and define a solution of \eqref{eq:1.1}.
Throughout this paper we often identify ${\bf R}^{N-1}$ with $\partial{\bf R}^N_+$. 
For any $x\in{\bf R}^N$ and $r>0$, let 
$$
B(x,r):=\{y\in{\bf R}^N\,:\,|x-y|<r\},
\quad
B_+(x,r):=\{(y',y_N)\in B(x,r)\,:\,y_N\ge 0\}.
$$
For any $L\ge 0$, we set 
\begin{equation*}
\begin{split}
 & D_L:=\{(x',x_N)\,:\,x'\in{\bf R}^{N-1},\,\, x_N\ge L^{1/2}\},\\
 & D_L':=\{(x',x_N)\,:\,x'\in{\bf R}^{N-1},\,\, 0\le x_N<L^{1/2}\}.
\end{split}
\end{equation*}
We remark that $D=D_0=\overline{{\bf R}^N_+}$. 
Let $\Gamma_N=\Gamma_N(x,t)$ be the Gauss kernel on ${\bf R}^N$, that is
\begin{equation}
\label{eq:1.4}
\Gamma_N(x,t):=(4\pi t)^{-\frac{N}{2}}\exp\left(-\frac{|x|^2}{4t}\right),
\quad 
x\in {\bf R}^N,\,\,\,t>0. 
\end{equation}
Let $G=G(x,y,t)$ be the Green function for the heat equation on ${\bf R}^N_+$ 
with the homogeneous Neumann boundary condition, that is 
\begin{equation}
\label{eq:1.5}
G(x,y,t):=\Gamma_N(x-y,t)+\Gamma_N(x-y_*,t),
\quad
x,y\in D,\,\,\,t>0,
\end{equation}
where $y_*=(y',-y_N)$ for $y=(y',y_N)\in D$. 
For any Radon measure $\mu$ in ${\bf R}^N$ with $\mbox{supp}\,\mu\subset D$, 
define
$$
[S(t)\mu](x):=\int_D G(x,y,t)\,d\mu(y),\quad x\in D,\,\,t>0. 
$$
For any locally integrable nonnegative function $\phi$ on $D$, 
we often identify $\phi$ with the Radon measure $\phi\,dx$.  
Then it follows that 
\begin{equation}
\label{eq:1.6}
\lim_{t\to+0}\|S(t)\eta-\eta\|_{L^\infty(D)}=0,
\qquad\eta\in C_0(D:[0,\infty)).
\end{equation}
\begin{definition}
\label{Definition:1.1}
Let $u$ be a nonnegative and continuous function in $D\times(0,T)$, 
where $0<T<\infty$. 
\newline
{\rm (i)} We say that $u$ is a solution of \eqref{eq:1.1} in $(0,T)$ 
if $u$ satisfies 
\begin{equation}
\label{eq:1.7}
u(x,t)=
\int_D G(x,y,\tau)u(y,t-\tau)\,dy
+\int_\tau^t\int_{{\bf R}^{N-1}}G(x,y',0,t-s)u(y',0,s)^p\,dy'\,ds
\end{equation}
for $(x,t)\in D\times (\tau,T)$ and $0<\tau<T$. 
\vspace{3pt}
\newline
{\rm (ii)} Let $\mu$ be a nonnegative measurable function in ${\bf R}^N_+$ 
or a Radon measure in ${\bf R}^N$ with $\mbox{supp}\,\mu\subset D$. 
We say that $u$ is a solution of \eqref{eq:1.1} and \eqref{eq:1.2} in $[0,T)$ 
if $u$ satisfies 
\begin{equation}
\label{eq:1.8}
u(x,t)=
\int_D G(x,y,t)\,d\mu+\int_0^t\int_{{\bf R}^{N-1}}G(x,y',0,t-s)u(y',0,s)^p\,dy'\,ds
\end{equation}
for $(x,t)\in D\times(0,T)$. 
If $u$ satisfies \eqref{eq:1.8} with $``="$ replaced by $``\ge"$, 
then $u$ is said to be a supersolution of \eqref{eq:1.1} and \eqref{eq:1.2} in $[0,T)$.  
\vspace{3pt}
\newline
{\rm (iii)} 
Let $u$ be a solution of \eqref{eq:1.1} and \eqref{eq:1.2} in $[0,T)$. 
We say that $u$ is a minimal solution of \eqref{eq:1.1} and \eqref{eq:1.2} in $[0,T)$ if
$u(x,t)\le v(x,t)$ in $D\times(0,T)$ 
for any solution $v$ of \eqref{eq:1.1} and \eqref{eq:1.2} in $[0,T)$. 
\end{definition}

Now we are ready to state our main results. 
In Theorem~\ref{Theorem:1.1} 
we show the existence and the uniqueness 
of the initial trace of the solution of \eqref{eq:1.1} 
and give necessary conditions on the initial trace.  
\begin{theorem}
\label{Theorem:1.1}
Let $u$ be a solution of \eqref{eq:1.1} in $(0,T)$, 
where $0<T<\infty$. 
Then there exists a unique Radon measure $\mu$ in ${\bf R}^N$ 
with $\mbox{supp}\,\mu\subset D$ such that 
\begin{equation}
\label{eq:1.9}
\lim_{t\to+0}
\int_D u(y,t)\phi(y)\,dy=\int_D \phi(y)\,d\mu(y),
\qquad
\phi\in C_0({\bf R}^N). 
\end{equation}
Furthermore, 
for any $\delta>0$, 
there exists $\gamma_1>0$ such that 
\begin{itemize}
  \item[{\rm (1)}] 
  $\displaystyle{\sup_{x\in D}\exp\left(-(1+\delta)\frac{x_N^2}{4T}\right)
  \mu(B(x,T^\frac{1}{2}))\le \gamma_1\,T^{\frac{N}{2}-\frac{1}{2(p-1)}}}$ 
  if $1<p<p_*$;
  \item[{\rm (2)}] 
  $\displaystyle{\sup_{x\in D}\exp\left(-(1+\delta)\frac{x_N^2}{4\sigma^2}\right)\mu(B(x,\sigma))
  \le \gamma_1\,
  \biggr[\log\biggr(e+\frac{T^{\frac{1}{2}}}{\sigma}\biggr)\biggr]^{-N}}$ 
  for  $0<\sigma<T^{\frac{1}{2}}$\\ if $p=p_*$; 
  \item[{\rm (3)}] 
  $\displaystyle{\sup_{x\in D}\exp\left(-(1+\delta)\frac{x_N^2}{4\sigma^2}\right)\mu(B(x,\sigma))
  \le \gamma_1\,
  \sigma^{N-\frac{1}{p-1}}}$ 
  for $0<\sigma<T^{\frac{1}{2}}$ if $p>p_*$.
\end{itemize} 
\end{theorem}
In Theorem~\ref{Theorem:1.2} 
we show that 
the initial trace of the solution of \eqref{eq:1.1} and \eqref{eq:1.2} coincides with its initial data. 
\begin{theorem}
\label{Theorem:1.2}
Let $\mu$ be a Radon measure in ${\bf R}^N$ with $\mbox{supp}\,\mu\subset D$. 
\begin{itemize}
  \item[{\rm (a)}] 
  Let $u$ be a solution of \eqref{eq:1.1} and \eqref{eq:1.2} in $[0,T)$ for some $T>0$. 
  Then \eqref{eq:1.9} holds. 
  \item[{\rm (b)}] 
  Let $u$ be a solution of \eqref{eq:1.1} in $(0,T)$ for some $T>0$. 
  Assume \eqref{eq:1.9}. 
  Then $u$ is a solution of \eqref{eq:1.1} and \eqref{eq:1.2} in $[0,T)$. 
\end{itemize}
\end{theorem}
Combining Theorem~\ref{Theorem:1.1} with Theorem~\ref{Theorem:1.2}, 
we obtain necessary conditions on the initial data for the solvability of problem~\eqref{eq:1.1} with \eqref{eq:1.2}. 
\begin{remark}
\label{Remark:1.1}
{\rm (i)} If $1<p\le p_*$ and $\mu\not\equiv 0$ on $D$, 
then problem~\eqref{eq:1.1} possesses no nonnegative global-in-time solutions. 
See {\rm\cite{DFL}} and {\rm\cite{GL}}. 
\vspace{3pt}
\newline
{\rm (ii)} 
Let $u$ be a solution of \eqref{eq:1.1} in $[0,\infty)$
and $1<p\le p_*$.  
It follows from assertions~{\rm (1)} and {\rm (2)} of Theorem~{\rm\ref{Theorem:1.1}} that 
the initial trace of $u$ must be identically zero in $D$. 
Then Theorem~{\rm\ref{Theorem:1.2}} 
leads the same conclusion as in Remark~{\rm\ref{Remark:1.1}}~{\rm (i)}. 
\end{remark}
Next we state our main results on sufficient conditions 
for the solvability of problem~\eqref{eq:1.1} with \eqref{eq:1.2}. 
In what follows, for any Radon measure in ${\bf R}^N$ and 
any bounded Borel set $E$, we denote by $|E|$ the Lebesgue measure of $E$ and  
set 
$$
\dashint_{E}\,d\mu=\frac{1}{|E|}\int_E\,d\mu. 
$$
\begin{theorem}
\label{Theorem:1.3}
Let $1<p<p_*$, $T>0$ and $\delta\in(0,1)$. 
Set $\lambda:=(1-\delta)/4T$. 
Then there exists $\gamma_2=\gamma_2(N,p,\delta)>0$ with the following property:
\begin{itemize}
 \item
 If $\mu$ is a Radon measure in ${\bf R}^N$ with $\mbox{supp}\,\mu\subset D$ satisfying
 \begin{equation}
 \label{eq:1.10}
 \sup_{x\in D}\,\dashint_{B(x,T^{\frac{1}{2}})}e^{-\lambda y_N^2}\,d\mu(y)
 \le\gamma_2T^{-\frac{1}{2(p-1)}},
 \end{equation}
 then there exists a solution $u$ of \eqref{eq:1.1} and \eqref{eq:1.2} in $[0,T)$ such that 
 $$
 0\le u(x,t)\le 2[S(t)\mu](x),\quad 
 (x,t)\in D\times(0,T). 
 $$
\end{itemize}
\end{theorem}
\begin{theorem}
\label{Theorem:1.4}
Let $1<\alpha<p$,  $T>0$ and $\delta\in(0,1)$. 
Set $\lambda:=(1-\delta)/4T$. 
Then there exists $\gamma_3=\gamma_3(N,p,\alpha,\delta)>0$ with the following property: 
\begin{itemize}
  \item Let $\mu_1$ be a Radon measure in ${\bf R}^N$ such that 
  $\mbox{supp}\,\mu_1\subset D_T$ and 
  \begin{equation}
  \label{eq:1.11}
  \sup_{x\in D_T}\,\dashint_{B(x,T^{\frac{1}{2}})}e^{-\lambda y_N^2}\,d\mu_1(y)
  \le\gamma_3T^{-\frac{1}{2(p-1)}}.
  \end{equation}
  Let $\mu_2$ be a nonnegative measurable function in ${\bf R}^N_+$ 
  such that $\mbox{supp}\,\mu_2\subset D_T'$ and 
  \begin{equation}
  \label{eq:1.12}
  \sup_{x\in D_T'}
  \left[\,\dashint_{B(x,\sigma)}
  \mu_2(y)^\alpha\,dy\,\right]^{\frac{1}{\alpha}}\le\gamma_3\sigma^{-\frac{1}{p-1}}
  \quad\mbox{for $0<\sigma<T^{\frac{1}{2}}$}.
  \end{equation}
  Then there exists a solution $u$ of \eqref{eq:1.1} and \eqref{eq:1.2} in $[0,T)$ 
  with $\mu=\mu_1+\mu_2$ such that 
  $$
  0\le u(x,t)\le 2[S(t)\mu_1](x)+2\left([S(t)\mu_2^\alpha](x)\right)^{\frac{1}{\alpha}},\quad 
  (x,t)\in D\times(0,T). 
  $$
\end{itemize}
\end{theorem}
\begin{theorem}
\label{Theorem:1.5}
Let $p=p_*$, $\beta>0$, $T>0$ and $\delta\in(0,1)$. Set $\lambda:=(1-\delta)/4T$ and 
\begin{equation}
\label{eq:1.13}
\Phi_\beta(s):=s[\log (e+s)]^\beta,
\quad
\rho(s):=
s^{-N}\biggr[\log\biggr(e+\frac{1}{s}\biggr)\biggr]^{-N}
\quad\mbox{for}\quad s>0. 
\end{equation}
Then there exists $\gamma_4=\gamma_4(N,\beta,\delta)>0$ with the following property:
\begin{itemize}
  \item Let $\mu_1$ be a Radon measure in ${\bf R}^N$ such that 
  $\mbox{supp}\,\mu_1\subset D_T$ and 
  \begin{equation}
  \label{eq:1.14}
  \sup_{x\in D_T}\,\dashint_{B(x,T^{\frac{1}{2}})}e^{-\lambda y_N^2}\,d\mu_1(y)
  \le\gamma_4T^{-\frac{1}{2(p-1)}}.
  \end{equation}
  Let $\mu_2$ be a nonnegative measurable function in ${\bf R}^N_+$ 
  such that $\mbox{supp}\,\mu_2\subset D_T'$ and 
  \begin{equation}
  \label{eq:1.15}
  \sup_{x\in D_T'}\Phi_\beta^{-1}\left[\,\dashint_{B(x,\sigma)}
  \Phi_\beta(T^\frac{1}{2(p-1)}\mu_2(y))\,dy\,\right]\le\gamma_4\rho(\sigma T^{-\frac{1}{2}})
  \quad\mbox{for $0<\sigma<T^{\frac{1}{2}}$}.
  \end{equation}
  Then there exists a solution of \eqref{eq:1.1} and \eqref{eq:1.2} in $[0,T)$ 
  with $\mu=\mu_1+\mu_2$ such that 
  $$
  0\le u(x,t)\le 2[S(t)\mu_1](x)+d\Phi_\beta^{-1}\left([S(t)\Phi_\beta(\mu_2)](x)\right),\quad 
  (x,t)\in D\times(0,T),
  $$
  where $d$ is a positive constant depending only on $p$ and $\beta$. 
\end{itemize} 
\end{theorem}
As a corollary of our theorems, we have: 
\begin{corollary}
\label{Corollary:1.1}
Let $\delta$ be the Delta function in ${\bf R}^N$ and $x_0\in D$. 
Let $\mu(y)=\delta(y-x_0)$ in ${\bf R}^N$. 
Then there exists a solution of \eqref{eq:1.1} and \eqref{eq:1.2} in $[0,T)$ for some $T>0$ if and only if, either 
$$
{\rm (i)}\quad x_0\in\partial{\bf R}^N_+\quad\mbox{and}\quad 1<p<p_*
\qquad\quad\mbox{or}\qquad\quad
{\rm (ii)}\quad x_0\in{\bf R}^N_+\quad\mbox{and}\quad p>1.
$$
\end{corollary}
See also Theorem~\ref{Theorem:6.3}. 
\vspace{3pt}

We develop the arguments in \cite{HI} and prove our theorems. 
Let $u$ be a solution of \eqref{eq:1.1} in $(0,T)$ for some $T>0$. 
By the same argument as in \cite{HI} we can prove the existence and the uniqueness of the initial trace 
of the solution~$u$. 
Furthermore, 
we study a lower estimate of the solution~$u$ near the boundary $\partial D$ 
by the use of $\|u(\tau)\|_{L^1(B_+(z,\rho))}$, where $z\in D$, $\rho\in(0,T^{1/2})$ and $\tau\in(0,T)$. 
(See Lemma~\ref{Lemma:3.1}.)
Combining this lower estimate with \cite[Lemma~2.1.2]{DFL}, 
we complete the proof of Theorem~\ref{Theorem:1.1} in the case $p\not=p_*$. 
For the case $p=p_*$, 
we obtain an integral inequality with respect to the quantity 
$$
\int_{\partial D}\Gamma_{N-1}(y',t)u(y',0,t)\,dy'
$$
(see \eqref{eq:3.15}). 
Then we apply a similar iteration argument as in \cite[Section~2]{LN} 
to obtain $\|u(\tau)\|_{L^1(B_+(z,\rho))}$, where $z\in D$, $\rho\in(0,T^{1/2})$ and $\tau\in(0,T)$. 
This completes the proof of Theorem~\ref{Theorem:1.1} in the case $p=p_*$. 
Theorem~\ref{Theorem:1.2} is proved by a similar argument as in the proof of \cite[Theorem~1.2]{HI} 
with the aid of Theorem~\ref{Theorem:1.1}. 
Furthermore, we prove a lemma 
on an estimate of an integral related to the nonlinear boundary condition 
(see Lemma~\ref{Lemma:5.1}) and apply the arguments in \cite{HI, IKS, RS} 
to prove Theorems~\ref{Theorem:1.3}--\ref{Theorem:1.5}. 
\vspace{5pt}

The rest of this paper is organized as follows. 
In Section~2 we recall some properties of the kernel $G=G(x,y,t)$ and 
prove some preliminary lemmas on the kernel $G$. 
In Section~3 we study the existence and the uniqueness of the initial trace.  
Furthermore, we obtain necessary conditions for the solvability of the solutions of \eqref{eq:1.1} and \eqref{eq:1.2}, 
and prove Theorem~\ref{Theorem:1.1}. 
In Section~4 we apply Theorem~\ref{Theorem:1.1} to prove Theorem~\ref{Theorem:1.2}. 
In Section~5 we obtain sufficient conditions on the initial data 
for the solvability of the solution of \eqref{eq:1.1} and \eqref{eq:1.2}, 
and prove Theorems~\ref{Theorem:1.3}, \ref{Theorem:1.4} and \ref{Theorem:1.5}. 
In Section~6, as an application of our theorems, 
we obtain some estimates of the life span of the solution of \eqref{eq:1.1} and \eqref{eq:1.2}. 
\section{Preliminaries}
In this section we recall some properties of the kernel $G=G(x,y,t)$ 
and prove preliminary lemmas.  
By \eqref{eq:1.5} we have
\begin{equation}
\label{eq:2.1}
\Gamma_N(x-y,t)\le G(x,y,t)\le 2\Gamma_N(x-y,t),
\qquad  x,y\in D,\,\,t>0. 
\end{equation}
It follows from \eqref{eq:1.4} and \eqref{eq:1.5} that 
\begin{equation}
\begin{split}
\label{eq:2.2}
G(x',x_N,y',0,t) & =G(y',0,x',x_N,t)=2\Gamma_N(x'-y',x_N,t)\\
 & =2(4\pi t)^{-\frac{1}{2}}\exp\left(-\frac{x_N^2}{4t}\right)\Gamma_{N-1}(x'-y',t)
\end{split}
\end{equation}
for $x\in D$, $y'\in{\bf R}^{N-1}$ and $t>0$. 
By the semigroup property of $S(t)$ we see that 
\begin{equation}
\label{eq:2.3}
\int_D G(x,y,t)G(y,z,s)\,dy=G(x,z,t+s)
\end{equation}
for $(x,t)$, $(z,s)\in D\times(0,\infty)$. 
Furthermore, 
we have the following two lemmas. 
In what follows, 
by the letter $C$
we denote generic positive constants 
and they may have different values also within the same line. 
\begin{lemma}
\label{Lemma:2.1}
Let $\mu$ be a Radon measure in ${\bf R}^N$ with $\mbox{supp}\,\mu\subset D$. 
If $[S(T)\mu](x)<\infty$ for some $x\in D$ and $T>0$, 
then $S(t)\mu$ is continuous in $D\times(0,T)$. 
\end{lemma} 
{\bf Proof.}
Assume $[S(T)\mu](x)<\infty$ for some $x\in D$ and $T>0$. Let $0<T'<T$. 
It follows from \eqref{eq:2.1} that 
\begin{equation*}
\begin{split}
\infty & >[S(T)\mu](x)\ge\int_D \Gamma_N(x-y,T)\,d\mu(y)\\
 & =(4\pi T)^{-\frac{N}{2}}\int_D
\exp\left(-\frac{|x-y|^2}{4T}\right)\,d\mu(y)
\ge C\int_D
\exp\left(-\frac{|y|^2}{4T'}\right)\,d\mu(y).
\end{split}
\end{equation*}
Then, applying the Lebesgue dominated convergence theorem, 
by \eqref{eq:2.1} 
we see that $S(t)\mu$ is continuous in $D\times(0,T')$. 
Since $T'$ is arbitrary,  the proof is complete. 
$\Box$
\begin{lemma}
\label{Lemma:2.2}
Let $T>0$ and $\mu$ be a Radon measure in ${\bf R}^N$. 
Let $\lambda\ge 0$ be such that $4\lambda T<1$. 
Assume that $\,\mbox{{\rm supp}}\,\mu\subset D_L$ for some $L\ge 0$. 
Then there exists $\gamma>0$ such that 
\begin{equation}
\label{eq:2.4}
[S(t)\mu](x)\le \gamma t^{-\frac{N}{2}}\exp\left(-\frac{L}{\gamma t}\right)
\sup_{z\in D_L}\int_{B(z,t^{\frac{1}{2}})}e^{-\lambda y_N^2}d\mu(y)
\end{equation}
for $x\in\partial{\bf R}^N_+$ and $0<t\le T$. 
\end{lemma}
{\bf Proof.} 
Let $x\in\partial{\bf R}^N_+$ and $0<t\le T$. 
Let $\lambda\ge 0$ be such that $4\lambda T<1$. 
By the Besicovitch covering lemma 
we can find an integer $m$ depending only on $N$ and 
a set $\{x_{k,i}\}_{k=1,\dots,m,\,i\in{\bf N}}\subset D_L$ such that 
\begin{equation}
\label{eq:2.5}
B_{k,i}\,\cap\,B_{k,j}=\emptyset\quad\mbox{if $i\not=j$},\qquad
D_L\subset\bigcup_{k=1}^m\bigcup_{i=1}^\infty B_{k,i},
\end{equation}
where $B_{k,i}:=\overline{B(x_{k,i},t^{\frac{1}{2}})}$. 
Then it follows from \eqref{eq:2.1} that 
\begin{equation}
\label{eq:2.6}
\begin{split}
 & [S(t)\mu](x)
 \le 2\sum_{k=1}^m\sum_{i=1}^\infty \int_{D_L\,\cap\, B_{k,i}}\Gamma_N(x-y,t)\,d\mu(y) \\
 & \le Ct^{-\frac{N}{2}}\sup_{k=1,\dots,m,\,i\in{\bf N}}\int_{B_{k,i}}e^{-\lambda y_N^2}d\mu(y)\,
 \sum_{k=1}^m \sum_{i=1}^\infty
 \sup_{y\in D_L\,\cap\, B_{k,i}}
\exp\left(-\frac{|x-y|^2}{4t}+\lambda y_N^2\right).
\end{split}
\end{equation}
On the other hand, 
for any $y\in D_L$ and $r>0$, 
there exists a set $\{y_\ell\}_{\ell=1}^{m'}\subset D_L$ such that 
\begin{equation}
\label{eq:2.7}
\overline{B(y,r)}\cap D_L\subset\bigcup_{\ell=1}^{m'}B(y_\ell,r)\cap D_L,
\end{equation}
where $m'$ is an integer depending only on $N$. 
This implies that 
\begin{equation}
\label{eq:2.8}
\int_{B_{k,i}}e^{-\lambda y_N^2}\,d\mu(y)
\le m'\sup_{z\in D_L}\int_{B(z,t^{\frac{1}{2}})}e^{-\lambda y_N^2}\,d\mu(y)
\end{equation}
for any $k\in\{1,\dots,m\}$, $i\in{\bf N}$ and $t\in(0,T]$. 
On the other hand, 
since $x\in\partial{\bf R}^N_+$ and 
\begin{equation*}
\begin{split}
|x'-y'|^2 & \ge (|x'-z'|-|z'-y'|)^2
=|x'-z'|^2-2|x'-z'||z'-y'|+|z'-y'|^2\\
 & \ge\frac{1}{2}|x'-z'|^2-|z'-y'|^2
\ge\frac{1}{2}|x'-z'|^2-4t,\\
y_N^2 & \ge (z_N-|y_N-z_N|)^2
=z_N^2-2z_N|y_N-z_N|+|y_N-z_N|^2\\
 & \ge\frac{1}{2}z_N^2-|y_N-z_N|^2
\ge\frac{1}{2}z_N^2-4t,
\end{split}
\end{equation*}
for $y$, $z\in B_{k,i}$,  
we have
\begin{equation*}
\begin{split}
\exp\left(-\frac{|x-y|^2}{4t}+\lambda y_N^2\right)
 & =\exp\left(-\frac{|x'-y'|^2}{4t}\right)\exp\left(-\frac{y_N^2}{4t}(1-4t\lambda)\right)\\
 & \le C\exp\left(-\frac{|x'-z'|^2}{8t}\right)
 \exp\left(-\delta\frac{z_N^2}{8t}\right)
\end{split}
\end{equation*}
for any $z\in B_{k,i}$, 
where $\delta:=1-4\lambda T>0$. 
This together with \eqref{eq:2.6} and \eqref{eq:2.8} implies that 
\begin{equation*}
\begin{split}
[S(t)\mu](x)
 & \le Ct^{-\frac{N}{2}}\sup_{z\in D_L}\int_{B(z,t^{\frac{1}{2}})}e^{-\lambda y_N^2}d\mu(y)\\
 & \qquad\qquad
\times\sum_{k=1}^m \sum_{i=1}^\infty\,
\dashint_{D_L\,\cap\, B_{k,i}}\exp\left(-\frac{|x'-z'|^2}{8t}\right)
\exp\left(-\delta\frac{z_N^2}{8t}\right)\,dz\\
 & \le Ct^{-N}\sup_{z\in D_L}\int_{B(z,t^{\frac{1}{2}})}e^{-\lambda y_N^2}d\mu(y)\,
 \int_{D_L} \exp\left(-\frac{|x'-z'|^2}{8t}\right)
 \exp\left(-\delta\frac{z_N^2}{8t}\right)\,dz\\
 & \le Ct^{-\frac{N}{2}}\exp\left(-\delta\frac{L}{16t}\right)
 \sup_{z\in D_L}\int_{B(z,t^{\frac{1}{2}})}e^{-\lambda y_N^2}d\mu(y).
\end{split}
\end{equation*}
Therefore we obtain \eqref{eq:2.4}. 
Thus Lemma~\ref{Lemma:2.2} follows.
$\Box$
\begin{lemma}
\label{Lemma:2.3}
Assume that there exists a supersolution $v$ of \eqref{eq:1.1} and \eqref{eq:1.2} 
in $[0,T)$ for some $T>0$. 
Then there exists a minimal solution of  \eqref{eq:1.1} and \eqref{eq:1.2} in $[0,T)$. 
\end{lemma}
{\bf Proof.} 
Since $v$ is a supersolution in $[0,T)$, we have 
$\infty>v(0,T')\ge [S(T')\mu](0)$
for any $T'\in(0,T)$. Then Lemma~\ref{Lemma:2.1} implies that 
$S(t)\mu\in C(D\times(0,T))$. 

Let $n\in\{1,2,\dots\}$. 
Set $u_{n,1}(x,t):=[S(t)\mu](x)$. 
Since $S(t)\mu\in C(D\times(0,T))$, 
we can define $u_{n,2}$ by 
$$
u_{n,2}(x,t):=[S(t)\mu](x)+\int_0^t\int_{{\bf R}^{N-1}}
G(x,y',0,t-s)\left(\min\{u_{n,1}(y',0,s),n\}\right)^p\,dy'\,ds 
$$
for $(x,t)\in D\times(0,T)$. 
Then it follows that 
$$
u_{n,2}\in C(D\times(0,T)),
\qquad
u_{n,2}(x,t)\le v(x,t)\quad\mbox{on}\quad D\times(0,T).
$$ 
By induction we define 
$u_{n,k}\in C(D\times(0,T))$
by 
$$
u_{n,k}(x,t):=[S(t)\mu](x)
+\int_0^t\int_{{\bf R}^{N-1}} G(x,y',0,t-s)\left(\min\{u_{n,k-1}(y',0,s),n\}\right)^p\,dy'
$$
for $(x,t)\in D\times(0,T)$, where $k=1,2,\dots$.
Furthermore, we see that 
\begin{equation*}
\begin{split}
u_{n,1}(x,t)\le u_{n,2}(x,t)\le\dots u_{n,k}(x,t)\le\dots\le v(x,t),\\
u_{1,k}(x,t)\le u_{2,k}(x,t)\le\dots u_{n,k}(x,t)\le\dots\le v(x,t),
\end{split}
\end{equation*}
for $(x,t)\in D\times(0,T)$. 
Then we deduce that 
the sequence $\{u_{n,k}\}$ is 
equibounded and equicontinuous with respect to $k$ and $n$ 
on any compact set $K\subset D\times(0,T)$ 
(see e.g., \cite[Section~6]{DB01} and \cite[Section~2]{IK}). 
By the Ascoli-Arzel\`a theorem and the diagonal argument 
we can find a function $u\in C(D\times(0,T))$ such that 
\begin{equation*}
u(x,t)=[S(t)\mu](x)
+\int_0^t\int_{{\bf R}^{N-1}} G(x,y',0,t-s)u(y',0,s)^p\,dy'\le v(x,t)
\end{equation*}
for $(x,t)\in D\times(0,T)$.  
This means that $u$ is a solution of \eqref{eq:1.1} and \eqref{eq:1.2} 
in $[0,T)$. Furthermore, we easily see that $u$ is a minimal solution of \eqref{eq:1.1} and \eqref{eq:1.2} 
in $[0,T)$. Thus Lemma~\ref{Lemma:2.3} follows. 
$\Box$
\vspace{5pt}

At the end of this section we state 
the following two lemmas on the initial trace of the solution of \eqref{eq:1.1}. 
These are proved by similar arguments as in the proofs of Lemmas~2.3 and 2.4 in \cite{HI}, respectively, 
and we left the proofs to the reader. 
\begin{lemma}
\label{Lemma:2.4}
Let $u$ be a solution of \eqref{eq:1.1} in $(0,T)$, where $0<T<\infty$. 
Then 
\begin{equation}
\label{eq:2.9}
\sup_{0<t<T-\epsilon}\,\int_{B_+(0,R)}u(y,t)\,dy<\infty
\end{equation}
for $R>0$ and $0<\epsilon<T$. 
Furthermore, 
there exists a unique Radon measure $\mu$ in ${\bf R}^N$ 
with $\mbox{{\rm supp}}\,\mu\subset D$
such that 
\begin{equation}
\label{eq:2.10}
\lim_{t\to +0}
\int_D u(y,t)\eta(y)\,dy=\int_D \eta(y)\,d\mu(y),
\quad \eta\in C_0({\bf R}^N).
\end{equation}
\end{lemma}
\begin{lemma}
\label{Lemma:2.5}
Let $\mu$ be a Radon measure in ${\bf R}^N$ with $\mbox{{\rm supp}}\,\mu\subset D$. 
Let $u$ be a solution of \eqref{eq:1.1} and \eqref{eq:1.2} in $[0,T)$ 
for some $0<T<\infty$. Then  \eqref{eq:2.10} holds for $\eta\in C_0({\bf R}^N)$. 
\end{lemma}
\section{Proof of Theorem~\ref{Theorem:1.1}}
In this section we prove Theorem~\ref{Theorem:1.1}. 
For this aim, we prepare the following lemma. 
\begin{lemma}
\label{Lemma:3.1}
Let $u$ be a solution of \eqref{eq:1.1} in $(0,T)$, 
where $0<T<\infty$. 
For any $\epsilon\in(0,1/2)$, 
there exists $\gamma_*>0$ such that 
\begin{equation}
\label{eq:3.1}
u(x+\overline{z},(1-\epsilon)T+\rho^2+\tau)
\ge \gamma_*\Gamma_N\left(x,\frac{T}{\gamma_*}\right)
\exp\left(-\frac{1+\epsilon}{1-\epsilon}\frac{z_N^2}{4T}\right)\int_{B_+(z,\rho)}u(y,\tau)\,dy
\end{equation}
for  $x$, $z\in D$, $\rho\in(0,(\epsilon T)^{1/2})$ and $\tau\in(0,(1-\epsilon)T)$, where $\overline{z}:=(z',0)$. 
Here the constant $\gamma_*$ depends only on $N$ and $\epsilon>0$. 
\end{lemma}
{\bf Proof.}
Let $z\in D$. 
We can assume, without loss of generality, that $z'=0$ and $\overline{z}=0$. 
Let $\epsilon\in(0,1/2)$, $0<\rho<(\epsilon T)^{1/2}$ and $\tau\in(0,(1-\epsilon)T)$. 
Since 
\begin{equation*}
\begin{split}
 & \min_{y\in B_+(z,\rho)}\Gamma_N(x-y,(1-\epsilon)T+\rho^2)
 \ge (4\pi((1-\epsilon)T+\rho^2))^{-\frac{N}{2}}\exp\left(-\frac{(|x|+|z|+\rho)^2}{4((1-\epsilon)T+\rho^2)}\right)\\
 & \ge (4\pi((1-\epsilon)T+\rho^2))^{-\frac{N}{2}}\exp\left(-\frac{C|x|^2+C\rho^2}{4((1-\epsilon)T+\rho^2)}\right)
 \exp\left(-\frac{(1+\epsilon)z_N^2}{4((1-\epsilon) T+\rho^2)}\right)\\
 & \ge C^{-1}\Gamma_N(x,CT)
 \exp\left(-\frac{1+\epsilon}{1-\epsilon}\frac{z_N^2}{4T}\right),
\end{split}
\end{equation*}
by \eqref{eq:1.7} and \eqref{eq:2.1} we obtain  
\begin{equation*}
\begin{split}
u(x,(1-\epsilon)T+\rho^2+\tau) & \ge\int_{B_+(z,\rho)}\Gamma_N(x-y,(1-\epsilon) T+\rho^2)u(y,\tau)\,dy\\
 & \ge C^{-1}\Gamma_N(x,CT)\exp\left(-\frac{1+\epsilon}{1-\epsilon}\frac{z_N^2}{4T}\right)
 \int_{B_+(z,\rho)}u(y,\tau)\,dy
\end{split}
\end{equation*}
for $x\in D$. 
This implies \eqref{eq:3.1}, 
and Lemma~\ref{Lemma:3.1} follows. 
$\Box$\vspace{5pt}
\newline
Next we recall the following lemma (see \cite[Lemma~2.1.2]{DFL}). 
\begin{lemma}
\label{Lemma:3.2}
Let $\mu\in C^1(D)$ be such that $\partial_{x_N}\mu\le 0$ in ${\bf R}^N_+$. 
Assume that there exists a solution of \eqref{eq:1.1} and \eqref{eq:1.2} 
in $[0,T)$ for some $T>0$. Then 
$$
[S(t)\mu](x',0)\le \gamma t^{-\frac{1}{2(p-1)}}
$$
holds for $x'\in{\bf R}^{N-1}$ and $t\in(0,T)$, where $\gamma$ is a constant depending only on $N$ and $p$. 
\end{lemma}

\noindent
Now we are ready to prove Theorem~\ref{Theorem:1.1}.
\vspace{5pt}
\newline
{\bf Proof of Theorem~\ref{Theorem:1.1}.}
By Lemma~\ref{Lemma:2.4} 
we can find a unique Radon measure $\mu$ in ${\bf R}^N$ with 
$\mbox{supp}\,\mu\subset D$ satisfying \eqref{eq:1.9}. 
So it suffices to prove assertions~(1), (2) and (3). 

Let $u$ be a solution of \eqref{eq:1.1} in $(0,T)$ for some $T>0$. 
Let $0<\sigma<T^{1/2}$ and $0<\epsilon<1/2$. 
Lemma~\ref{Lemma:3.1} implies that 
\begin{equation}
\label{eq:3.2}
u(x+\overline{z},(1-\epsilon)\sigma^2+\rho^2+\tau)
\ge \gamma_*\Gamma_N\left(x,\frac{\sigma^2}{\gamma_*}\right)
\exp\left(-\frac{1+\epsilon}{1-\epsilon}\frac{z_N^2}{4\sigma^2}\right)\int_{B_+(z,\rho)}u(y,\tau)\,dy
\end{equation}
for  $x$, $z\in D$, $\rho\in(0,\epsilon^{1/2}\sigma)$ and $\tau\in(0,(1-\epsilon)\sigma^2)$, 
where $\gamma_*$ is as in Lemma~\ref{Lemma:3.1}. 
\vspace{3pt}
\newline
{\bf Proof of assertions~(1) and (3).}
Since $\tilde{u}(x,t):=u(x+\overline{z},t+(1-\epsilon)\sigma^2+\rho^2+\tau)$ is a solution of \eqref{eq:1.1} 
in $(0,\epsilon\sigma^2-\rho^2-\tau)$,  
by Lemma~\ref{Lemma:2.3} and \eqref{eq:3.2} 
we can find a minimal solution $w$ of \eqref{eq:1.1} in $[0,\epsilon\sigma^2-\rho^2-\tau)$ 
with 
$$
w(x,0)=\gamma_*\Gamma_N\left(x,\frac{\sigma^2}{\gamma_*}\right)
\exp\left(-\frac{1+\epsilon}{1-\epsilon}\frac{z_N^2}{4\sigma^2}\right)\int_{B_+(z,\rho)}u(y,\tau)\,dy,
\qquad 
x\in D. 
$$
Then it follows from Lemma~\ref{Lemma:3.2} that 
\begin{equation}
\label{eq:3.3}
\begin{split}
 & Ct^{-\frac{1}{2(p-1)}}\ge [S(t)w(0)](0)\\
 & =\gamma_*\Gamma_N\left(0,t+\frac{\sigma^2}{\gamma_*}\right)
\exp\left(-\frac{1+\epsilon}{1-\epsilon}\frac{z_N^2}{4\sigma^2}\right)\int_{B_+(z,\rho)}u(y,\tau)\,dy\\
 & =\gamma_*(4\pi)^{-\frac{N}{2}}\left(t+\frac{\sigma^2}{\gamma_*}\right)^{-\frac{N}{2}}
 \exp\left(-\frac{1+\epsilon}{1-\epsilon}\frac{z_N^2}{4\sigma^2}\right)\int_{B_+(z,\rho)}u(y,\tau)\,dy
\end{split}
\end{equation}
for $0<t<\epsilon\sigma^2-\rho^2-\tau$. 

Let $0<\rho'<\rho$. Let $\zeta\in C_0({\bf R}^N)$ be such that 
$$
\zeta=1\quad\mbox{on}\quad B(z,\rho'),
\qquad
0\le \zeta\le 1\quad{in}\quad{\bf R}^N,
\qquad
\zeta=0\quad\mbox{outside}\quad B(z,\rho). 
$$
By Lemma~\ref{Lemma:2.4} we have 
\begin{equation}
\label{eq:3.4}
\limsup_{\tau\to +0}\int_{B_+(z,\rho)}u(y,\tau)\,dy
\ge\limsup_{\tau\to +0}\int_D u(y,\tau)\zeta\,dy
 =\int_D \zeta\,d\mu(y)
\ge\int_{B_+(z,\rho')}\,d\mu(y). 
\end{equation}
Since $\rho'$ is arbitrary, 
by \eqref{eq:3.3} and \eqref{eq:3.4} we obtain 
\begin{equation}
\label{eq:3.5}
\gamma_*(4\pi)^{-\frac{N}{2}}\left(t+\frac{\sigma^2}{\gamma_*}\right)^{-\frac{N}{2}}
\exp\left(-\frac{1+\epsilon}{1-\epsilon}\frac{z_N^2}{4\sigma^2}\right)\int_{B_+(z,\rho)}\,d\mu\le
Ct^{-\frac{1}{2(p-1)}} 
\end{equation}
for $z\in D$, $\rho\in(0,\epsilon^{1/2}\sigma)$ and $0<t<\epsilon\sigma^2-\rho^2$. 
Setting $\rho=(\epsilon/2)^{1/2}\sigma$ and $t=\epsilon\sigma^2/4$, 
we obtain  
\begin{equation}
\label{eq:3.6}
\exp\left(-\frac{1+\epsilon}{1-\epsilon}\frac{z_N^2}{4\sigma^2}\right)
\int_{B_+(z,(\epsilon/2)^{1/2}\sigma)}\,d\mu
\le C\sigma^{N-\frac{1}{p-1}}
\end{equation}
for $z\in D$ and $\sigma\in(0,T^{1/2})$. 

On the other hand, for any $z\in D$, 
we can find $\{z_\ell\}_{\ell=1}^{m'}\subset D$ such that 
\begin{equation}
\label{eq:3.7}
B_+(z,\sigma)\subset \bigcup_{\ell=1}^{m'}B_+(z_\ell,(\epsilon/2)^{1/2}\sigma). 
\end{equation}
Here $m'$ is independent of $z$. 
We can assume, without loss of generality, that
$B_+(z,\sigma)\cap B_+(z_\ell,(\epsilon/2)^{1/2}\sigma)\not=\emptyset$. 
Then $(z_\ell)_N\le z_N+2\sigma$ and it follows that
\begin{equation*}
\begin{split}
 & \exp\left(-\frac{(1+\epsilon)^2}{1-\epsilon}\frac{z_N^2}{4\sigma^2}
+\frac{1+\epsilon}{1-\epsilon}\frac{(z_\ell)_N^2}{4\sigma^2}\right)\\
 & \le
\exp\left(-\frac{(1+\epsilon)^2}{1-\epsilon}\frac{z_N^2}{4\sigma^2}
+\frac{1+\epsilon}{1-\epsilon}\frac{(1+\epsilon)z_N^2+C(2\sigma)^2}{4\sigma^2}\right)
\le\exp\left(C\frac{1+\epsilon}{1-\epsilon}\right)\le C.
\end{split}
\end{equation*}
This together with \eqref{eq:3.6} and \eqref{eq:3.7} implies that 
\begin{equation}
\label{eq:3.8}
\begin{split}
 & \exp\left(-\frac{(1+\epsilon)^2}{1-\epsilon}\frac{z_N^2}{4\sigma^2}\right)
\int_{B_+(z,\sigma)}\,d\mu\\
 & \le\sum_{\ell=1}^{m'}
\exp\left(-\frac{(1+\epsilon)^2}{1-\epsilon}\frac{z_N^2}{4\sigma^2}
+\frac{1+\epsilon}{1-\epsilon}\frac{(z_\ell)_N^2}{4\sigma^2}\right)\\
 & \qquad
 \times
 \exp\left(-\frac{1+\epsilon}{1-\epsilon}\frac{(z_\ell)_N^2}{\sigma^2}\right)
 \int_{B_+(z_\ell,(\epsilon/2)^{1/2}\sigma)}\,d\mu
 \le C\sigma^{N-\frac{1}{p-1}}
\end{split}
\end{equation}
for $z\in D$ and $0<\sigma<T^{1/2}$. 

Let $\delta>0$. Taking a sufficiently small $\epsilon\in(0,1/2)$ if necessary, 
we have $(1+\epsilon)^2/(1-\epsilon)\le 1+\delta$.
Then \eqref{eq:3.8} implies assertions~(1) and (3). 
\vspace{5pt}
\newline
\noindent
{\bf Proof of assertion~(2).} 
Let $p=p_*$. 
Set $\rho\in(0,(\epsilon/2)^{1/2}\sigma)$ and $\tau\in(0,(1-\epsilon)\sigma^{1/2})$. 
For any $z=(z',z_N)\in D$, 
set  
$$
v(x,t):=u(x+\overline{z},t+(1-\epsilon)\sigma^2+\rho^2)
$$
for $x\in D$ and $t\in(0,T-(1-\epsilon)\sigma^2-\rho^2)$, where $\overline{z}=(z',0)$. 
Since $v$ is a solution of \eqref{eq:1.1} in $(0,T-(1-\epsilon)\sigma^2-\rho^2)$, 
we have
\begin{equation}
\label{eq:3.9}
v(x,t)=\int_D G(x,y,t-\tau)v(y,\tau)\,dy
+\int_\tau^t\int_{{\bf R}^{N-1}}G(x,y',0,t-s)v(y',0,s)^p\,dy'\,ds
\end{equation}
for $x\in D$ and $0<\tau<t<T-(1-\epsilon)\sigma^2-\rho^2$. 
In particular, 
for any $0<T'<(T-(1-\epsilon)\sigma^2-\rho^2)/2$, by \eqref{eq:2.2} we have 
\begin{equation*}
\begin{split}
\infty>v(0,2T')
 & \ge\int_\tau^{T'}\int_{{\bf R}^{N-1}}G(0,y',0,2T'-s)v(y',0,s)^p\,dy'\,ds\\
 & =2\int_\tau^{T'}(4\pi(2T'-s))^{-\frac{1}{2}}
 \int_{{\bf R}^{N-1}}\Gamma_{N-1}(y',2T'-s)v(y',0,s)^p\,dy'\,ds\\
 & \ge 2\int_\tau^{T'}(4\pi(2T'-s))^{-\frac{N}{2}}(4 \pi s)^{\frac{N-1}{2}}
 \int_{{\bf R}^{N-1}}\Gamma_{N-1}(y',s)v(y',0,s)^p\,dy'\,ds
\end{split}
\end{equation*}
for $0<\tau<T'$.  
This together with the Jensen inequality that 
\begin{equation*}
\begin{split}
\infty & >2\int_\tau^{T'}(4\pi(2T'-s))^{-\frac{N}{2}}(4 \pi s)^{\frac{N-1}{2}}
\left(\int_{{\bf R}^{N-1}}\Gamma_{N-1}(y',s)v(y',0,s)\,dy'\right)^p\,ds
\end{split}
\end{equation*}
for $0<\tau<T'$. 
Since $T'$ is arbitrary, we see that 
\begin{equation}
\label{eq:3.10}
V(t):=\int_{{\bf R}^{N-1}}\Gamma_{N-1}(y',t)v(y',0,t)\,dy'<\infty
\end{equation}
for almost all $t\in(0,(T-(1-\epsilon)\sigma^2-\rho^2)/2)$. 

It follows from \eqref{eq:3.2} that 
\begin{equation}
\label{eq:3.11}
\begin{split}
 & \int_D G(x,y,t-\tau)v(y,\tau)\,dy
 =\int_D G(x,y,t-\tau)u(y+\overline{z},(1-\epsilon)\sigma^2+\rho^2+\tau)\,dy\\
 & \ge \gamma_*\exp\left(-\frac{1+\epsilon}{1-\epsilon}\frac{z_N^2}{4\sigma^2}\right)\int_{B_+(z,\rho)}u(y,\tau)\,dy\,
\int_D G(x,y,t-\tau)\Gamma_N\left(y,\frac{\sigma^2}{\gamma_*}\right)\,dy\\
 & =\gamma_*\exp\left(-\frac{1+\epsilon}{1-\epsilon}\frac{z_N^2}{4\sigma^2}\right)\int_{B_+(z,\rho)}u(y,\tau)\,dy\cdot
\Gamma_N\left(x,t-\tau+\frac{\sigma^2}{\gamma_*}\right)
\end{split}
\end{equation}
for $x\in D$ and $0<\tau<t<T-(1-\epsilon)\sigma^2-\rho^2$, 
where $\gamma_*$ is as in Lemma~\ref{Lemma:3.1}. 
Setting
$$
M_\tau:=\gamma_*\exp\left(-\frac{1+\epsilon}{1-\epsilon}\frac{z_N^2}{4\sigma^2}\right)\int_{B_+(z,\rho)}u(y,\tau)\,dy, 
$$
by \eqref{eq:3.9}, \eqref{eq:3.10} and \eqref{eq:3.11} we obtain 
\begin{equation}
\label{eq:3.12}
\begin{split}
\infty>V(t) & \ge M_\tau\int_{{\bf R}^{N-1}}\Gamma_{N-1}(x',t)\Gamma_N\left(x',0,t-\tau+\frac{\sigma^2}{\gamma_*}\right)\,dx'\\
 & +\int_{{\bf R}^{N-1}}\int_\tau^t\int_{{\bf R}^{N-1}}G(x',0,y',0,t-s)\Gamma_{N-1}(x',t)v(y',0,s)^p\,dy'\,ds\,dx'
\end{split}
\end{equation}
for almost all $0<\tau<t<(T-(1-\epsilon)\sigma^2-\rho^2)/2$. 
It follows from $0<\rho^2<\epsilon\sigma^2/2$ that 
\begin{equation}
\label{eq:3.13}
\begin{split}
 & \int_{{\bf R}^{N-1}}\Gamma_{N-1}(x',t)\Gamma_N\left(x',0,t-\tau+\frac{\sigma^2}{\gamma_*}\right)\,dx'\\
 & =\int_{{\bf R}^{N-1}}\Gamma_{N-1}(x',t)\biggr(4\pi\biggr(t-\tau+\frac{\sigma^2}{\gamma_*}\biggr)\biggr)^{-\frac{1}{2}}
 \Gamma_{N-1}\biggr(x',t-\tau+\frac{\sigma^2}{\gamma_*}\biggr)\,dx'\\
 & =\biggr(4\pi\biggr(t-\tau+\frac{\sigma^2}{\gamma_*}\biggr)\biggr)^{-\frac{1}{2}}
 \Gamma_{N-1}\biggr(0,2t-\tau+\frac{\sigma^2}{\gamma_*}\biggr)\\
 & =\biggr(4\pi\biggr(t-\tau+\frac{\sigma^2}{\gamma_*}\biggr)\biggr)^{-\frac{1}{2}}
 \biggr(4\pi\biggr(2t-\tau+\frac{\sigma^2}{\gamma_*}\biggr)\biggr)^{-\frac{N-1}{2}}
 \ge ct^{-\frac{N}{2}}
\end{split}
\end{equation}
for $0<\tau<\epsilon\sigma^2/4<t<(T-(1-\epsilon)\sigma^2-\rho^2)/2$, 
where $c$ is a positive constant depending only on $N$ and $\epsilon$. 
Furthermore, 
by \eqref{eq:2.2} and the Jensen inequality we have 
\begin{equation}
\label{eq:3.14}
\begin{split}
 & \int_{{\bf R}^{N-1}}\int_\tau^t\int_{{\bf R}^{N-1}}G(x',0,y',0,t-s)\Gamma_{N-1}(x',t)v(y',0,s)^p\,dy'\,ds\,dx'\\
 & =\int_\tau^t\int_{{\bf R}^{N-1}}
 2[4\pi (t-s)]^{-\frac{1}{2}}\\
 & \qquad\qquad
 \times
 \biggr[\int_{{\bf R}^{N-1}}\Gamma_{N-1}(x',t)\Gamma_{N-1}(x'-y',t-s)\,dx'\biggr]v(y',0,s)^p\,dy'\,ds\\
 & =\int_\tau^t\int_{{\bf R}^{N-1}}
 2[4\pi (t-s)]^{-\frac{1}{2}}\Gamma_{N-1}(y',2t-s)v(y',0,s)^p\,dy'\,ds\\
 & \ge\int_\tau^t\int_{{\bf R}^{N-1}}
 2[4\pi (t-s)]^{-\frac{1}{2}}\biggr(\frac{s}{2t}\biggr)^{\frac{N-1}{2}}\Gamma_{N-1}(y',s)v(y',0,s)^p\,dy'\,ds\\
 & \ge\int_\tau^t
 2[4\pi (t-s)]^{-\frac{1}{2}}\biggr(\frac{s}{2t}\biggr)^{\frac{N-1}{2}}V(s)^p\,ds\\
 & =2^{-\frac{N-1}{2}}\pi^{-\frac{1}{2}}t^{-\frac{N-1}{2}}
 \int_\tau^t
 (t-s)^{-\frac{1}{2}}s^{\frac{N-1}{2}}V(s)^p\,ds. 
\end{split}
\end{equation} 
Therefore, 
by \eqref{eq:3.12}, \eqref{eq:3.13} and \eqref{eq:3.14} 
we obtain 
\begin{equation}
\label{eq:3.15}
V(t)\ge cM_\tau t^{-\frac{N}{2}}
+2^{-\frac{N-1}{2}}\pi^{-\frac{1}{2}}t^{-\frac{N-1}{2}}
\int_{\epsilon\sigma^2/4}^t (t-s)^{-\frac{1}{2}}s^{\frac{N-1}{2}}V(s)^p\,ds
\end{equation}
for $0<\tau<\epsilon\sigma^2/4$ 
and almost all $t\in(\epsilon\sigma^2/4,(T-(1-\epsilon)\sigma^2-\rho^2)/2)$. 

Set $a_1:=c$ and $\omega_1(t):=a_1M_\tau t^{-\frac{N}{2}}$.
Define
$$
\omega_{n+1}(t):=2^{-\frac{N-1}{2}}\pi^{-\frac{1}{2}}t^{-\frac{N-1}{2}}
\int_{\epsilon\sigma^2/4}^t (t-s)^{-\frac{1}{2}}s^{\frac{N-1}{2}}\omega_n(s)^p\,ds,
\qquad n=1,2,\dots,
$$
for $t>\epsilon\sigma^2/4$.  
Then it follows that 
\begin{equation}
\label{eq:3.16}
\infty>V(t)\ge\omega_{n+1}(t)\ge a_{n+1}M_\tau^{p^n}
t^{-\frac{N}{2}}\biggr[\log\left(\frac{4t}{\epsilon\sigma^2}\right)\biggr]^{\frac{p^n-1}{p-1}}
\end{equation}
for almost all $t\in(\epsilon\sigma^2/4,(T-(1-\epsilon)\sigma^2-\rho^2)/2)$ and $n=0,1,2,\dots$. 
Here $\{a_n\}$ is a sequence defined by 
\begin{equation}
\label{eq:3.17}
a_{n+1}:=2^{-\frac{N-1}{2}}\pi^{-\frac{1}{2}}a_n^p\frac{p-1}{p^n-1},
\quad n=1,2,\dots. 
\end{equation}
Indeed, \eqref{eq:3.16} holds with $n=0$. 
Furthermore, if \eqref{eq:3.16} holds for some $n\in\{0,1,2,\dots\}$, then, 
by \eqref{eq:3.15} we have 
\begin{equation*}
\begin{split}
 & \infty>V(t)\ge\omega_{n+2}(t)=2^{-\frac{N-1}{2}}\pi^{-\frac{1}{2}}t^{-\frac{N-1}{2}}
\int_{\epsilon\sigma^2/4}^t (t-s)^{-\frac{1}{2}}s^{\frac{N-1}{2}}\omega_{n+1}(s)^p\,ds\\
 & 
 =2^{-\frac{N-1}{2}}\pi^{-\frac{1}{2}}a_{n+1}^pM_\tau^{p^{n+1}}t^{-\frac{N-1}{2}}
 \int_{\epsilon\sigma^2/4}^t (t-s)^{-\frac{1}{2}}s^{\frac{N-1}{2}}s^{-\frac{N}{2}\left(1+\frac{1}{N}\right)}
 \biggr[\log\left(\frac{4s}{\epsilon\sigma^2}\right)\biggr]^{\frac{p^{n+1}-p}{p-1}}\,ds\\
 & 
 \ge 2^{-\frac{N-1}{2}}\pi^{-\frac{1}{2}}a_{n+1}^pM_\tau^{p^{n+1}}t^{-\frac{N}{2}}
 \int_{\epsilon\sigma^2/4}^t s^{-1} \biggr[\log\left(\frac{4s}{\epsilon\sigma^2}\right)\biggr]^{\frac{p^{n+1}-p}{p-1}}\,ds\\
 & 
 =2^{-\frac{N-1}{2}}\pi^{-\frac{1}{2}}a_{n+1}^pM_\tau^{p^{n+1}}t^{-\frac{N}{2}}
 \frac{p-1}{p^{n+1}-1}\biggr[\log\left(\frac{4t}{\epsilon\sigma^2}\right)\biggr]^{\frac{p^{n+1}-1}{p-1}}\\
  &
 =a_{n+2}M_\tau^{p^{n+1}}t^{-\frac{N}{2}}\biggr[\log\left(\frac{4t}{\epsilon\sigma^2}\right)\biggr]^{\frac{p^{n+1}-1}{p-1}}
\end{split}
\end{equation*}
for almost all $t\in(\epsilon\sigma^2/4,(T-(1-\epsilon)\sigma^2-\rho^2)/2)$. 
This means that \eqref{eq:3.16} holds for $n+1$. 
Thus \eqref{eq:3.16} holds for all $n\in\{0,1,2,\dots\}$. 

On the other hand, 
similarly to \cite[Lemma~2.20~(i)]{LN} 
(see also (3.26) in \cite{HI}), 
we can find $b>0$ such that 
$$
a_n\ge b^{p^n},
\qquad n=1,2,\dots. 
$$
This together with \eqref{eq:3.16} implies that 
\begin{equation*}
\begin{split}
\infty>V(t)\ge\omega_{n+1}(t)
 & \ge b^{p^{n+1}}M_\tau^{p^n}t^{-\frac{N}{2}}\biggr[\log\left(\frac{4t}{\epsilon\sigma^2}\right)\biggr]^{\frac{p^n-1}{p-1}}\\
 & =t^{-\frac{N}{2}}\biggr[\log\left(\frac{4t}{\epsilon\sigma^2}\right)\biggr]^{-\frac{1}{p-1}}
 \left(b^pM_\tau\biggr[\log\left(\frac{4t}{\epsilon\sigma^2}\right)\biggr]^{\frac{1}{p-1}}\right)^{p^n}
\end{split}
\end{equation*}
for almost all $t\in(\epsilon\sigma^2/4,(T-(1-\epsilon)\sigma^2-\rho^2)/2)$ and $n=1,2,\dots$. 
Then it follows that 
$$
M_\tau\le b^{-p}\biggr[\log\left(\frac{4t}{\epsilon\sigma^2}\right)\biggr]^{-\frac{1}{p-1}}
=b^{-p}\biggr[\log\left(\frac{4t}{\epsilon\sigma^2}\right)\biggr]^{-N},
$$
which implies that 
$$
\exp\left(-\frac{1+\epsilon}{1-\epsilon}\frac{z_N^2}{4\sigma^2}\right)\int_{B_+(z,\rho)}u(y,\tau)\,dy
\le (b^p\gamma_*)^{-1}\biggr[\log\left(\frac{4t}{\epsilon\sigma^2}\right)\biggr]^{-N},
\,\,\,\,
z\in D,\,\,\tau\in\left(0,\frac{\epsilon\sigma^2}{4}\right),
$$
for $t\in(\epsilon\sigma^2/4,(T-(1-\epsilon)\sigma^2-\rho^2)/2)$.  
Then, similarly to \eqref{eq:3.4}, we obtain 
\begin{equation}
\label{eq:3.18}
\exp\left(-\frac{1+\epsilon}{1-\epsilon}\frac{z_N^2}{4\sigma^2}\right)\int_{B_+(z,\rho)}\,d\mu(y)
\le (b^p\gamma_*)^{-1}\biggr[\log\left(\frac{4t}{\epsilon\sigma^2}\right)\biggr]^{-N}
\end{equation}
for $z\in D$ and $t\in(\epsilon\sigma^2/4,(T-(1-\epsilon)\sigma^2-\rho^2)/2)$. 

Set $\rho=(\epsilon/4)^{1/2}\sigma$.
Consider the case where $0<\sigma^2\le T/2$. It follows that 
$$
\frac{T-(1-\epsilon)\sigma^2-\rho^2}{2}>\frac{T-\sigma^2}{2}\ge\frac{T}{4}. 
$$
Setting $t=T/4$, by \eqref{eq:3.18} we have 
\begin{equation}
\label{eq:3.19}
\begin{split}
 & \exp\left(-\frac{1+\epsilon}{1-\epsilon}\frac{z_N^2}{4\sigma^2}\right)\int_{B_+(z,(\epsilon/4)^{1/2}\sigma)}\,d\mu(y)\\
 & \le (b^p\gamma_*)^{-1}\biggr[\log\left(\frac{T}{\epsilon\sigma^2}\right)\biggr]^{-N}
 \le C\biggr[\log\left(e+\frac{T^{\frac{1}{2}}}{\sigma}\right)\biggr]^{-N},
\quad z\in D. 
\end{split}
\end{equation}
On the other hand, in the case where $T/2<\sigma^2<T$, 
we have 
$$
\frac{T-(1-\epsilon)\sigma^2-\rho^2}{2}
\ge\frac{\epsilon\sigma^2-\rho^2}{2}
=\frac{3}{8}\epsilon\sigma^2,
\qquad
1<\frac{T}{\sigma^2}<2.
$$ 
Then, taking a sufficiently small $\epsilon\in(0,1/2)$ if necessary, 
we set
$$
t=\frac{5}{16}\epsilon\sigma^2\in\left(\frac{\epsilon\sigma^2}{4},\frac{(T-(1-\epsilon)\sigma^2-\rho^2)}{2}\right)
$$ 
and by \eqref{eq:3.18} we obtain 
\begin{equation}
\label{eq:3.20}
\begin{split}
 & \exp\left(-\frac{1+\epsilon}{1-\epsilon}\frac{z_N^2}{4\sigma^2}\right)\int_{B_+(z,(\epsilon/4)^{1/2}\sigma)}\,d\mu(y)\\
 & \le (b^p\gamma_*)^{-1}\biggr[\log\left(\frac{5}{4}\right)\biggr]^{-N}
\le C\biggr[\log\left(e+\frac{T^{\frac{1}{2}}}{\sigma}\right)\biggr]^{-N}, 
 \quad z\in D. 
\end{split}
\end{equation}
Combining \eqref{eq:3.19} and \eqref{eq:3.20} and applying the same argument as in \eqref{eq:3.8}, 
we obtain 
$$
\exp\left(-\frac{(1+\epsilon)^2}{1-\epsilon}\frac{z_N^2}{4\sigma^2}\right)\mu(B(z,\sigma))
\le C\biggr[\log\left(e+\frac{T^{\frac{1}{2}}}{\sigma}\right)\biggr]^{-N},
\quad
z\in D,\,\,\sigma\in(0,T^{1/2}).
$$
Finally, similarly to the proof of assertions~(1) and (3), 
for any $\delta>0$, we take a sufficiently small $\epsilon\in(0,1/2)$ to obtain assertion~(2). 
Thus Theorem~\ref{Theorem:1.1} follows.
$\Box$ 
\section{Proof of Theorem~\ref{Theorem:1.2}}
We modify the proof of \cite[Theorem~1.2]{HI} 
to prove Theorem~\ref{Theorem:1.2}. 
\vspace{5pt}
\newline
{\bf Proof of Theorem~\ref{Theorem:1.2}.}
By Lemma~\ref{Lemma:2.5} it suffices to prove Theorem~\ref{Theorem:1.2}~(b). 
Let $u$ be a solution of \eqref{eq:1.1} in $(0,T)$, where $0<T<\infty$.  
By \eqref{eq:3.3} there exists $\gamma>0$ such that 
\begin{equation}
\label{eq:4.1}
\exp\left(-4\gamma\frac{z_N^2}{T}\right)\int_{B_+(z,(T/4)^{\frac{1}{2}})} u(y,\tau)\,dy
\le CT^{\frac{N}{2}-\frac{1}{2(p-1)}}
\end{equation}
for $z\in D$ and $\tau\in(0,T/8)$. 
Let $t\in(0,T)$. 
For any $n=1,2,\dots$,  
by the Besicovitch covering lemma we can find an integer $m$
depending only on $N$ and 
a set $\{x_{k,i}\}_{k=1,\dots,m,\,i\in{\bf N}}\subset D\setminus B(0,nt^{1/2})$ such that 
\begin{equation}
\label{eq:4.2}
B_{k,i}\cap B_{k,j}=\emptyset\quad\mbox{if $i\not=j$}
\qquad\mbox{and}\qquad
D\setminus B(0,nt^\frac{1}{2})\subset\bigcup_{k=1}^m\bigcup_{i=1}^\infty B_{k,i},
\end{equation}
where $B_{k,i}:=\overline{B_+(x_{k,i},t^{1/2})}$. 
Since 
\begin{equation*}
\begin{split}
 & \sup_{y\in B_{k,i}} \exp\left(4\gamma\frac{(x_{k,i})_N^2}{T}\right)G(y,t-\tau)\\
 & \le 2(4\pi(t-\tau))^{-\frac{N}{2}}\sup_{y\in B_{k,i}}
\exp\left(4\gamma\frac{[z_N+|(x_{k,i})_N-z_N|]^2}{T}\right)
\exp\left(-\frac{|y|^2}{4(t-\tau)}\right)\\
 & \le Ct^{-\frac{N}{2}}
\exp\left(8\gamma\frac{z_N^2}{T}\right)
\sup_{y\in B_{k,i}}
\exp\left(-\frac{(|z|-|z-y|)^2}{4t}\right)\\
& \le Ct^{-\frac{N}{2}}\exp\left(8\gamma\frac{z_N^2}{T}\right)\exp\left(-\frac{|z|^2}{8t}\right)
\end{split}
\end{equation*}
for $z\in B_{k,i}$ and $0<\tau<t/2$, 
by \eqref{eq:4.1} and \eqref{eq:4.2} we obtain  
\begin{equation*}
\begin{split}
 & \sup_{0<\tau<t/2}
 \int_{D\setminus B(0,nt^\frac{1}{2})}G(y,t-\tau)u(y,\tau)\,dy
 \le\sum^m_{k=1}\sum^\infty_{i=1}\,\sup_{0<\tau<t/2}\int_{B_{k,i}}G(y,t-\tau)u(y,\tau)\,dy\\
 & \qquad\quad
 \le C\,\sup_{0<\tau<t/2}\sup_{z\in D}\,
 \exp\left(-4\gamma\frac{z_N^2}{T}\right)
 \int_{B_+(z,(T/4)^\frac{1}{2})}u(y,\tau)\,dy\\
 & \qquad\qquad\qquad
 \times\sum^m_{k=1}\sum^\infty_{i=1}\sup_{0<\tau<t/2}\sup_{y\in B_{k,i}}
 \exp\left(4\gamma\frac{(x_{k,i})_N^2}{T}\right)G(y,t-\tau)\\
  & \qquad\quad
  \le CT^{\frac{N}{2}-\frac{1}{2(p-1)}}t^{-\frac{N}{2}}
 \sum^m_{k=1}\sum^\infty_{i=1}
 \dashint_{B_{k,i}}
 \exp\left(8\gamma\frac{z_N^2}{T}\right)\exp\left(-\frac{|z|^2}{8t}\right)\,dz
\end{split}
\end{equation*}
for $0<t<T/4$ and $0<\tau<t/2$. 
Then, taking a sufficiently small $t>0$, we see that 
$$
\sup_{0<\tau<t/2}
\int_{D\setminus B(0,nt^\frac{1}{2})}G(y,t-\tau)u(y,\tau)\,dy\\
\le Ct^{-N}\sum^m_{k=1}\sum^\infty_{i=1}
\int_{B_{k,i}}\exp\left(-\frac{|z|^2}{16t}\right)\,dz,
$$
which together with \eqref{eq:4.2} implies that 
\begin{equation}
\label{eq:4.3}
\begin{split}
  & \sup_{0<\tau<t/2}
 \int_{D\setminus B(0,nt^\frac{1}{2})}G(y,t-\tau)u(y,\tau)\,dy\\
  & 
  \le Ct^{-N}
 \int_{D\setminus B(0,(n-1)t^{\frac{1}{2}})}\exp\left(-\frac{|z|^2}{16t}\right)\,dz\to 0
\end{split}
\end{equation}
as $n\to\infty$. 
Similarly, by using Theorem~\ref{Theorem:1.1}, instead of \eqref{eq:3.3}, 
we see that 
\begin{equation}
\label{eq:4.4}
\lim_{n\to\infty}\sup_{0<\tau<t/2}\int_{D\setminus B(0,nt^\frac{1}{2})}G(y,t-\tau)\,d\mu=0
\end{equation}
for all sufficiently small $t>0$. 

Let $\eta_n\in C_0({\bf R}^N)$ be such that 
$$
0\le\eta_n\le 1\quad\mbox{in}\quad{\bf R}^N,
\qquad
\eta_n=1\quad\mbox{on}\quad B(0,nt^\frac{1}{2}),
\qquad
\eta_n=0\quad\mbox{outside}\quad B(0,2nt^\frac{1}{2}). 
$$
Then we have 
\begin{equation}
\label{eq:4.5}
\begin{split}
 & \left|\int_D G(y,t-\tau)u(y,\tau)\,dy-\int_D G(y,t)\,d\mu(y)\right|\\
 & \le\left|\int_D G(y,t)u(y,\tau)\eta_n(y)\,dy - \int_D G(y,t)\eta_n(y)\,d\mu(y)\right|\\
 & \qquad
 +\left|\int_D[G(y,t-\tau)-G(y,t)]u(y,\tau)\eta_n(y)\,dy\right|\\
 & \qquad\qquad
 +\int_{D\setminus B(0,nt^\frac{1}{2})}G(y,t-\tau)u(y,\tau)\,dy
 +\int_{D\setminus B(0,nt^\frac{1}{2})} G(y,t)\,d\mu(y)
\end{split}
\end{equation}
for $n=1,2,\dots$ and $\tau\in(0,t/2)$. 
By Lemma~\ref{Lemma:2.4} we see that 
\begin{equation}
\label{eq:4.6}
\lim_{\tau\to+0}\,
\left[\int_D G(y,t)u(y,\tau)\eta_n(y)\,dy - \int_D G(y,t)\eta_n(y)\,d\mu(y)\right]=0.
\end{equation}
Furthermore, by Lemma~\ref{Lemma:2.4} we have 
\begin{equation}
\label{eq:4.7}
\begin{split}
 & \lim_{\tau\to+0}\,\left|\int_D [G(y,t-\tau)-G(y,t)]u(y,\tau)\eta_n(y)\,dy\right|\\
 & 
\le \sup_{y\in B(0,2nt^{\frac{1}{2}}),s\in(t/2,t)}\,|\partial_t G(y,s)|
\,\limsup_{\tau\to+0}\,
\biggr[\tau\int_{B_+(0,2nt^\frac{1}{2})}u(y,\tau)\,dy\biggr]=0.
\end{split}
\end{equation}
By \eqref{eq:4.5}, \eqref{eq:4.6} and \eqref{eq:4.7} we see that 
\begin{equation*}
\begin{split}
 & \limsup_{\tau\to+0}\,
 \left|\int_D G(y,t-\tau)u(y,\tau)\,dy - \int_D G(y,t)\,d\mu(y)\right|\\
 & \le
\sup_{0<\tau<t/2}\,\int_{D\setminus B(0,nt^\frac{1}{2})}G(y,t-\tau)u(y,\tau)\,dy
 +\int_{D\setminus B(0,nt^\frac{1}{2})} G(y,t)\,d\mu(y)
\end{split}
\end{equation*}
for $n=1,2,\dots$. 
This together with \eqref{eq:4.3} and \eqref{eq:4.4} implies that 
$$
\lim_{\tau\to+0}\,
\left|\int_D G(y,t-\tau)u(y,\tau)\,dy - \int_D G(y,t)\,d\mu(y)\right|=0.
$$
This together with Definition~\ref{Definition:1.1}~(i) implies that 
$u$ is a solution of \eqref{eq:1.1} and \eqref{eq:1.2} in $[0,T)$. 
Thus Theorem~\ref{Theorem:1.2}~(b) follows, and the proof is complete. 
$\Box$
\section{Proof of Theorems~\ref{Theorem:1.3}, \ref{Theorem:1.4} and \ref{Theorem:1.5}}
We prove Theorems~\ref{Theorem:1.3}, \ref{Theorem:1.4} and \ref{Theorem:1.5} 
by modifying the arguments in \cite{HI, IKS, RS}. 
Lemma~\ref{Lemma:5.1} is a key lemma in our proofs.
\begin{lemma}
\label{Lemma:5.1}
Let $\Psi$ be a nonnegative convex function on $[0,\infty)$. 
Let $\mu$ be a nonnegative measurable function in $D$ such that 
$[S(T)\Psi(\mu)](x_0)<\infty$
for some $x_0\in D$ and $T>0$. 
Define
\begin{equation*}
\begin{split}
W(x,t):=& \Psi^{-1}\left([S(t)\Psi(\mu)](x)\right),\qquad
w(x',t):=W(x',0,t),\\
F(x,t):= & \int_0^t\int_{{\bf R}^{N-1}}G(x,y',0,t-s)
w(y',s)^p\,dy'\,ds,
\end{split}
\end{equation*}
for $x'\in{\bf R}^{N-1}$, $x\in D$ and $t\in(0,T)$. 
Then there exists $\gamma>0$ such that
\begin{equation}
\label{eq:5.1}
\begin{split}
F(x,t)
 & \le \gamma t^{\frac{1}{2}}
 \left\|\frac{\Psi(W(t))}{W(t)}\right\|_{L^\infty(D)}
 W(x,t)\\
 & \qquad
\times\int_0^t (t-s)^{-\frac{1}{2}}s^{-\frac{1}{2}}
\left\|\frac{w(s)^p}{\Psi(w(s))}\right\|_{L^\infty({\bf R}^{N-1})}\,ds\\
\end{split}
\end{equation}
for $x\in D$ and $0<t<T$. 
\end{lemma}
{\bf Proof.}
Since $[S(T)\Psi(\mu)](x_0)<\infty$ for some $x_0\in D$ and $T>0$, 
by Lemma~\ref{Lemma:2.1} we can define $w=w(x',t)$ and $F=F(x,t)$ 
for $x'\in{\bf R}^{N-1}$, $x\in D$ and $t\in(0,T)$. 
It follows that 
\begin{equation}
\label{eq:5.2}
F(x,t)\le\int_0^t
\left\|\frac{w(s)^p}{\Psi(w(s))}\right\|_{L^\infty({\bf R}^{N-1})}
\int_{{\bf R}^{N-1}}G(x,y',0,t-s)
[S(s)\Psi(\mu)](y',0)\,dy'\,ds. 
\end{equation}
On the other hand, by \eqref{eq:2.2} we have  
\begin{equation*}
\begin{split}
 & \int_{{\bf R}^{N-1}}G(x,y',0,t-s)
[S(s)\Psi(\mu)](y',0)\,dy'\\
 & =\int_{{\bf R}^{N-1}}G(x',x_N,y',0,t-s)
\int_D G(y',0,z',z_N,s)\Psi(\mu(z))\,dz\,dy'\\
 & =\int_D\int_{{\bf R}^{N-1}}
 2(4\pi(t-s))^{-\frac{1}{2}}\Gamma_{N-1}(x'-y',t-s)\exp\left(-\frac{x_N^2}{4(t-s)}\right)\\
 & \qquad\qquad
 \times
 2(4\pi s)^{-\frac{1}{2}}\exp\left(-\frac{z_N^2}{4s}\right)\Gamma_{N-1}(y'-z',s)
 \Psi(\mu(z',z_N))\,dy'\,dz
\end{split}
\end{equation*} 
for $x\in{\bf R}^{N-1}$ and $0<s<t<T$. 
Then we have 
\begin{equation*}
\begin{split}
 & \int_{{\bf R}^{N-1}}G(x,y',0,t-s)
[S(s)\Psi(\mu)](y',0)\,dy'\\
 & \le\exp\left(-\frac{x_N^2}{4t}\right)\int_D
 2(4\pi(t-s))^{-\frac{1}{2}}\Gamma_{N-1}(x'-z',t)\\
 & \qquad\qquad
 \times
 2(4\pi s)^{-\frac{1}{2}}\exp\left(-\frac{z_N^2}{4s}\right)
 \Psi(\mu(z',z_N))\,dz\\
 & \le\exp\left(-\frac{x_N^2}{4t}\right)
 \int_D2(4\pi t)^{-\frac{1}{2}}\exp\left(-\frac{z_N^2}{4t}\right)
 \Gamma_{N-1}(x'-z',t)\\
 & \qquad\qquad
 \times (t-s)^{-\frac{1}{2}}t^{\frac{1}{2}}2(4\pi s)^{-\frac{1}{2}}\Psi(\mu(z',z_N))\,dz\\
 & =\exp\left(-\frac{x_N^2}{4t}\right)(t-s)^{-\frac{1}{2}}t^{\frac{1}{2}}2(4\pi s)^{-\frac{1}{2}}
 \int_D G(x',0,z',z_N,t)\Psi(\mu(z))\,dz\\
 & =\exp\left(-\frac{x_N^2}{4t}\right)\pi^{-\frac{1}{2}}t^{\frac{1}{2}}
 (t-s)^{-\frac{1}{2}}s^{-\frac{1}{2}}[S(t)\Psi(\mu)](x',0)
\end{split}
\end{equation*}
for $x\in{\bf R}^{N-1}$ and $0<s<t<T$. 
This together with \eqref{eq:5.2} implies that 
\begin{equation}
\label{eq:5.3}
\begin{split}
F(x,t) & \le Ct^{\frac{1}{2}}\exp\left(-\frac{x_N^2}{4t}\right)[S(t)\Psi(\mu)](x',0)\\
 & \qquad\quad
 \times\int_0^t \left\|\frac{w(s)^p}{\Psi(w(s))}\right\|_{L^\infty({\bf R}^{N-1})}
(t-s)^{-\frac{1}{2}}s^{-\frac{1}{2}}\,ds
\end{split}
\end{equation}
for $x\in D$ and $t\in(0,T)$. 
Furthermore, it follows from \eqref{eq:2.1} and \eqref{eq:2.2} that 
\begin{equation*}
\begin{split}
 & \exp\left(-\frac{x_N^2}{4t}\right)[S(t)\Psi(\mu)](x',0)\\
 & =2(4\pi t)^{-\frac{N}{2}}\int_D 
\exp\left(-\frac{|x'-y'|^2}{4t}-\frac{x_N^2+y_N^2}{4t}\right)\Psi(\mu)\,dy\\
 & \le 2(4\pi t)^{-\frac{N}{2}}\int_D 
\exp\left(-\frac{|x'-y'|^2}{4t}-\frac{|x_N-y_N|^2}{4t}\right)\Psi(\mu)\,dy\\
& =2\int_D \Gamma_N(x-y,t)\Psi(\mu)\,dy
\le 2\int_D G(x,y,t)\Psi(\mu)\,dy
\le 2[S(t)\Psi(\mu)](x)
\end{split}
\end{equation*}
for $x=(x',x_N)\in D$ and $t\in(0,T)$. 
This together with \eqref{eq:5.3} implies that 
\begin{equation*}
\begin{split}
 F(x,t)
  & \le Ct^{\frac{1}{2}}[S(t)\Psi(\mu)](x)
 \int_0^t \left\|\frac{w(s)^p}{\Psi(w(s))}\right\|_{L^\infty({\bf R}^{N-1})}
 (t-s)^{-\frac{1}{2}}s^{-\frac{1}{2}}\,ds\\
 & \le Ct^{\frac{1}{2}}
 \left\|\frac{\Psi(W(t))}{W(t)}\right\|_{L^\infty(D)}
 W(x,t)\int_0^t \left\|\frac{w(s)^p}{\Psi(w(s))}\right\|_{L^\infty({\bf R}^{N-1})}
 (t-s)^{-\frac{1}{2}}s^{-\frac{1}{2}}\,ds
 \end{split}
\end{equation*}
for $x\in D$ and $t\in(0,T)$. 
Thus we obtain \eqref{eq:5.1}, and the proof is complete.
$\Box$
\vspace{5pt}
\newline
{\bf Proof of Theorem~\ref{Theorem:1.3}.}
It suffices to consider the case $T=1$. 
Indeed, for any solution~$u$ of \eqref{eq:1.1} in $[0,T)$, where $0<T<\infty$, 
$T^{1/2(p-1)}u(T^{1/2} x,T t)$ is a solution of \eqref{eq:1.1} in $[0,1)$. 

Let $\delta\in(0,1)$ and set $\lambda:=(1-\delta)/4$. 
Assume \eqref{eq:1.10}. Then $[S(1)\mu](0)<\infty$ 
and Lemma~\ref{Lemma:2.1} implies that $S(t)\mu\in C(D\times(0,1))$. 
We define
$$
\overline{u}(x,t):=2[S(t)\mu](x),
\quad
W(x,t):=[S(t)\mu](x),
\quad
w(x',t):=W(x',0,t),
$$
for $x\in D$, $x'\in{\bf R}^{N-1}$ and $t\in(0,1)$. 
It follows from Lemma~\ref{Lemma:2.2} with $L=0$ and \eqref{eq:1.10} that  
\begin{equation}
\label{eq:5.4}
\|w(t)\|_{L^\infty({\bf R}^{N-1})}
\le C t^{-\frac{N}{2}}
\sup_{z\in D}\int_{B(z,t^{\frac{1}{2}})}e^{-\lambda y_N^2}\,d\mu(y)\\
\le C\gamma_2t^{-\frac{N}{2}}
\end{equation}
for $t\in(0,1)$.
Applying Lemma~\ref{Lemma:5.1} with $\Psi(\tau)=\tau$, 
by \eqref{eq:5.4} we obtain 
\begin{equation}
\label{eq:5.5}
\begin{split}
 & \int_0^t\int_{{\bf R}^{N-1}}G(x,y',0,t-s)w(y',s)^p\,dy'\,ds\\
 & \le Ct^{\frac{1}{2}}
 W(x,t)
\int_0^t (t-s)^{-\frac{1}{2}}s^{-\frac{1}{2}}
\|w(s)\|^{p-1}_{L^\infty({\bf R}^{N-1})}\,ds\\
 & \le C\gamma_2^{p-1}t^{\frac{1}{2}} W(x,t)
 \int_0^t (t-s)^{-\frac{1}{2}}s^{-\frac{1}{2}-\frac{N}{2}(p-1)}\,ds
 \le C\gamma_2^{p-1}W(x,t)
\end{split}
\end{equation}
for $x=(x',x_N)\in D$ and $t\in(0,1)$. 
Taking a sufficiently small $\gamma_2>0$ if necessary, 
by \eqref{eq:5.5} we obtain 
\begin{equation*}
\begin{split}
 & \int_D G(x,y,t)\,d\mu(y)+\int_0^t\int_{{\bf R}^{N-1}}G(x,y',0,t-s)\overline{u}(y',0,s)^p\,dy'\,ds\\
 & =[S(t)\mu](x)+2^p\int_0^t\int_{{\bf R}^{N-1}}G(x,y',0,t-s)w(y',s)^p\,dy'\,ds\\
 & \le [S(t)\mu](x)+W(x,t)=2[S(t)\mu](x)=\overline{u}(x,t)
\end{split}
\end{equation*}
for $(x,t)\in D\times(0,1)$. 
This means that $\overline{u}$ is a supersolution of \eqref{eq:1.1} and \eqref{eq:1.2} in $[0,1)$. 
Therefore, by Lemma~\ref{Lemma:2.3} 
we can find a solution of \eqref{eq:1.1} and \eqref{eq:1.2}  in $[0,1)$ 
such that 
$0\le u(x,t)\le\overline{u}(x,t)=2[S(t)\mu](x)$ in $D\times(0,1)$. 
Thus Theorem~\ref{Theorem:1.3} follows. 
$\Box$
\vspace{5pt}
\newline
{\bf Proof of Theorem~\ref{Theorem:1.4}.}
Similarly to the proof of Theorem~\ref{Theorem:1.3}, 
it suffices to consider the case $T=1$. 
Let $\delta\in(0,1)$ and set $\lambda:=(1-\delta)/4$.
Assume \eqref{eq:1.11} and \eqref{eq:1.12}. 
Then Lemma~\ref{Lemma:2.1} implies that $S(t)\mu_1$, $S(t)\mu_2^\alpha\in C(D\times(0,1))$ and we define
\begin{equation*}
\begin{split}
 & \overline{u}(x,t):=
2[S(t)\mu_1](x)+2\left([S(t)\mu_2^\alpha](x)\right)^{\frac{1}{\alpha}},\\
 & W_1(x,t):=[S(t)\mu_1](x),\qquad
 W_2(x,t):=\left([S(t)\mu_2^\alpha](x)\right)^{\frac{1}{\alpha}},\\
 & w_1(x',t):=W_1(x',0,t),
\qquad
w_2(x',t):=W_2(x',0,t),
\end{split}
\end{equation*}
for $x\in D$, $x'\in{\bf R}^{N-1}$ and $t\in(0,1)$. 
Then it follows from the Jensen inequality that 
\begin{equation}
\label{eq:5.6}
\begin{split}
 & \int_D G(x,y,t)\,d\mu(y)+\int_0^t\int_{{\bf R}^{N-1}}G(x,y',0,t-s)\overline{u}(y',0,s)^p\,dy'\,ds\\
 & \le [S(t)\mu_1](x)+2^{2p-1}\int_0^t\int_{{\bf R}^{N-1}}G(x,y',0,t-s)w_1(y',s)^p\,dy'\,ds\\
 & \quad
 +\left([S(t)\mu_2^\alpha](x)\right)^{\frac{1}{\alpha}}
 +2^{2p-1}\int_0^t\int_{{\bf R}^{N-1}}G(x,y',0,t-s)w_2(y',s)^p\,dy'\,ds
\end{split}
\end{equation}
for $(x,t)\in D\times(0,1)$. 

On the other hand, 
by Lemma~\ref{Lemma:2.2} and \eqref{eq:1.11} we see that 
\begin{equation}
\label{eq:5.7}
\begin{split}
\|w_1(t)\|_{L^\infty({\bf R}^{N-1})} & \le C t^{-\frac{N}{2}}\exp\left(-\frac{1}{Ct}\right)
\sup_{z\in D_1}\int_{B(z,t^{\frac{1}{2}})}e^{-\lambda y_N^2}\,d\mu_1(y)\\
 & \le C\gamma_3 t^{-\frac{N}{2}}\exp\left(-\frac{1}{Ct}\right)
 \le C\gamma_3,
 \qquad t\in(0,1].
\end{split}
\end{equation}
Applying Lemma~\ref{Lemma:5.1} with $\Psi(\tau)=\tau$, 
by \eqref{eq:5.7} we obtain 
\begin{equation}
\label{eq:5.8}
\begin{split}
 & \int_0^t\int_{{\bf R}^{N-1}}G(x,y',0,t-s)w_1(y',s)^p\,dy'\,ds\\
 & \le Ct^{\frac{1}{2}}W_1(x,t)
\int_0^t (t-s)^{-\frac{1}{2}}s^{-\frac{1}{2}}
\|w_1(s)\|^{p-1}_{L^\infty({\bf R}^{N-1})}\,ds
\le C\gamma_3^{p-1}W_1(x,t)
\end{split}
\end{equation}
for $(x,t)\in D\times(0,1)$. 
Similarly, 
by Lemma~\ref{Lemma:2.2} and \eqref{eq:1.12} with $\sigma=t^{1/2}$ 
we see that 
\begin{equation}
\label{eq:5.9}
\begin{split}
\|W_2(t)\|_{L^\infty(D)}
\le \left[C\sup_{z\in D_1}\dashint_{B(z,t^{\frac{1}{2}})}e^{-\lambda y_N^2}\mu_2^\alpha(y)\,dy\right]^{\frac{1}{\alpha}}
\le C\gamma_3 t^{-\frac{1}{2(p-1)}}
\end{split}
\end{equation}
for $t\in(0,1)$.
Applying Lemma~\ref{Lemma:5.1} with $\Psi(\tau)=\tau^\alpha$, 
by \eqref{eq:5.9} we obtain 
\begin{equation}
\label{eq:5.10}
\begin{split}
 & \int_0^t\int_{{\bf R}^{N-1}}G(x,y',0,t-s)w_2(y',s)^p\,dy'\,ds\\
 & \le Ct^{\frac{1}{2}}\|W_2(t)\|_{L^\infty(D)}^{\alpha-1}W_2(x,t)
 \int_0^t (t-s)^{-\frac{1}{2}}s^{-\frac{1}{2}}
 \|w_2(s)\|_{L^\infty({\bf R}^{N-1})}^{p-\alpha}\,ds\\
 & \le C\gamma_3^{p-1}t^{\frac{1}{2}}t^{-\frac{\alpha-1}{2(p-1)}}W_2(x,t)
 \int_0^t (t-s)^{-\frac{1}{2}}s^{-\frac{1}{2}}s^{-\frac{p-\alpha}{2(p-1)}}\,ds
 \le C\gamma_3^{p-1}W_2(x,t)
\end{split}
\end{equation}
for $(x,t)\in D\times(0,1)$. 
Taking a sufficiently small $\gamma_3>0$ if necessary, 
by \eqref{eq:5.6}, \eqref{eq:5.8} and \eqref{eq:5.10} we obtain 
\begin{equation}
\label{eq:5.11}
\begin{split}
 & \int_D G(x,y,t)\,d\mu(y)+\int_0^t\int_{{\bf R}^{N-1}}G(x,y',0,t-s)\overline{u}(y',0,s)^p\,dy'\,ds\\
 & \le [S(t)\mu_1](x)+W_1(x,t)
 +\left([S(t)\mu_2^\alpha](x)\right)^{\frac{1}{\alpha}}
 +W_2(x,t)=\overline{u}(x,t)
\end{split}
\end{equation}
for $(x,t)\in D\times(0,1)$. 
This means that $\overline{u}$ is a supersolution of \eqref{eq:1.1} and \eqref{eq:1.2} in $[0,1)$. 
Therefore, by Lemma~\ref{Lemma:2.3} we can find a solution of \eqref{eq:1.1} and \eqref{eq:1.2} in $[0,1)$ 
such that 
$0\le u(x,t)\le\overline{u}(x,t)$ in $D\times(0,1)$.
Thus Theorem~\ref{Theorem:1.4} follows. 
$\Box$
\vspace{5pt}
\newline
{\bf Proof of Theorem~\ref{Theorem:1.5}.}
It suffices to consider the case $T=1$. 
Let $\Phi_\beta(s)$ and $\rho(s)$ be as in \eqref{eq:1.13}. 
Let $h\ge e$ be such that  
\begin{itemize}
  \item[{\rm (a)}] 
  $\Phi(\tau):=\tau[\log (h+\tau)]^\beta$ is positive and convex in $(0,\infty)$;
  \item[{\rm (b)}] 
  $\tau^p/\Phi(\tau)$ and $\Phi(\tau)/\tau$ are monotone increasing in $(0,\infty)$.
\end{itemize}
Let $\delta\in(0,1)$ and set $\lambda:=(1-\delta)/4$.
Assume \eqref{eq:1.14} and \eqref{eq:1.15}. 
Then Lemma~\ref{Lemma:2.1} implies that $S(t)\mu_1$, $S(t)\Phi(\mu_2)\in C(D\times(0,1))$. 
We define
\begin{equation*}
\begin{split}
 & \overline{u}(x,t):=
2[S(t)\mu_1](x)+d\Phi_\beta^{-1}\left([S(t)\Phi_\beta(\mu_2)](x)\right),\\
 & W_1(x,t):=[S(t)\mu_1](x),
\quad
 W_2(x,t):=\Phi^{-1}\left([S(t)\Phi(\mu_2)](x)\right),\\
 & w_1(x',t):=W_1(x',0,t),
\quad\,\,\,
 w_2(x',t):=W_2(x',0,t),
\end{split}
\end{equation*}
for $x\in D$, $x'\in{\bf R}^{N-1}$ and $t\in(0,1)$. 
Here $d$ is a positive constant to be chosen later. 
Since 
$$
\overline{u}(x,t)\le 2[S(t)\mu_1](x)+cdW_2(x,t)\quad\mbox{in}\quad D\times(0,1)
$$
for some $c>0$, 
similarly to \eqref{eq:5.6}, it follows from the Jensen inequality that
\begin{equation}
\label{eq:5.12}
\begin{split}
 & \int_{{\bf R}^N_+} G(x,y,t)\,d\mu(y)+\int_0^t\int_{{\bf R}^{N-1}}G(x,y',0,t-s)\overline{u}(y',0,s)^p\,dy'\,ds\\
 & \le [S(t)\mu_1](x)+2^{2p-1}\int_0^t\int_{{\bf R}^{N-1}}G(x,y',0,t-s)w_1(y',s)^p\,dy'\,ds\\
 & \quad
 +W_2(x,t)+2^{p-1}(cd)^p\int_0^t\int_{{\bf R}^{N-1}}G(x,y',0,t-s)w_2(y',s)^p\,dy'\,ds
\end{split}
\end{equation}
for $(x,t)\in D\times(0,1)$. 
Furthermore, similarly to \eqref{eq:5.8}, 
we obtain 
\begin{equation}
\label{eq:5.13}
\int_0^t\int_{{\bf R}^{N-1}}G(x,y',0,t-s)w_1(y',s)^p\,dy'\,ds
\le C\gamma_4^{p-1}W_1(x,t)
\end{equation}
for $(x,t)\in D\times(0,1)$. 

On the other hand, 
it follows from Lemma~\ref{Lemma:2.2} that 
\begin{equation}
\label{eq:5.14}
\begin{split}
 & \|\Phi(W_2(t))\|_{L^\infty(D)}
=\|S(t)\Phi(\mu_2)\|_{L^\infty(D)}\\
 & \qquad\quad
 \le Ct^{-\frac{N}{2}}
\sup_{z\in D}\int_{B(z,t^{\frac{1}{2}})}
\Phi(\mu_2)\,dy
\le Ct^{-\frac{N}{2}}
\sup_{z\in D}\int_{B(z,t^{\frac{1}{2}})}
\Phi_\beta(\mu_2)\,dy\\
 & \qquad\quad
 \le C\Phi_\beta\left(\gamma_4\rho(t^{\frac{1}{2}})\right)
 \le C\Phi\left(\gamma_4\rho(t^{\frac{1}{2}})\right)
\end{split}
\end{equation}
for $t\in(0,1)$. 
This together with property~(b) of $\Phi$ implies that 
\begin{equation}
\label{eq:5.15}
\left\|\frac{w_2(s)^p}{\Phi(w_2(s))}\right\|_{L^\infty({\bf R}^{N-1})}
\le\frac{\|w_2(s)\|_{L^\infty({\bf R}^{N-1})}^p}{\Phi(\|w_2(s)\|_{L^\infty({\bf R}^{N-1})})}
\le\frac{[{\Phi^{-1}}(C\Phi(\gamma_4\rho(s^\frac{1}{2})))]^p}
{C\Phi(\gamma_4\rho(s^\frac{1}{2}))}
\end{equation}
for $s\in(0,1)$.
Furthermore, we have
\begin{equation}
\label{eq:5.16}
\Phi(\gamma_4\rho(s^\frac{1}{2}))
=\gamma_4\rho(s^\frac{1}{2})[\log(h+\gamma_4\rho(s^\frac{1}{2}))]^\beta
\left\{
\begin{array}{l}
\ge \displaystyle{C^{-1}\gamma_4 s^{-\frac{N}{2}}\biggr[\log\biggr(e+\frac{1}{s}\biggr)\biggr]^{-N+\beta}},\vspace{7pt}\\
\le \displaystyle{C\gamma_4 s^{-\frac{N}{2}}\biggr[\log\biggr(e+\frac{1}{s}\biggr)\biggr]^{-N+\beta}},
\end{array}
\right.\end{equation}
for $s>0$.
Since 
$\Phi^{-1}(\tau)\le C\tau[\log(e+\tau)]^{-\beta}$ for $\tau>0$, 
it follows that 
\begin{equation}
\label{eq:5.17}
\Phi^{-1}(C\Phi(\gamma\rho(s^\frac{1}{2})))\le C\gamma_4 s^{-\frac{N}{2}}
\biggr[\log\biggr(e+\frac{1}{s}\biggr)\biggr]^{-N}
\end{equation}
for $s\in(0,1)$. 
By \eqref{eq:5.15}, \eqref{eq:5.16} and \eqref{eq:5.17} 
we obtain  
\begin{equation}
\label{eq:5.18}
\begin{split}
\left\|\frac{w_2(s)^p}{\Phi(w_2(s))}\right\|_{L^\infty({\bf R}^N)}
\le C\gamma_4^{\frac{1}{N}}s^{-\frac{1}{2}}\biggr[\log\biggr(e+\frac{1}{s}\biggr)\biggr]^{-1-\beta}
\end{split}
\end{equation}
for $s\in(0,1)$. Similarly, by \eqref{eq:5.14} and property~(b) we have
\begin{equation}
\label{eq:5.19}
\left\|\frac{\Phi(W_2(t))}{W_2(t)}\right\|_{L^\infty(D)}
\le\frac{C\Phi(\gamma_4\rho(t^\frac{1}{2}))}
{\Phi^{-1}(C\Phi(\gamma_4\rho(t^\frac{1}{2})))}
\le C\biggr[\log\biggr(e+\frac{1}{t}\biggr)\biggr]^\beta
\end{equation}
for $t\in(0,1)$.
By \eqref{eq:5.18} and \eqref{eq:5.19} we apply Lemma~\ref{Lemma:5.1} with $\Psi(\tau)=\Phi(\tau)$ 
to obtain 
\begin{equation}
\label{eq:5.20}
\begin{split}
 & \int_0^t\int_{{\bf R}^{N-1}}G(x,y',0,t-s)
w_2(y',s)^p\,dy'\,ds\\
 & \le C\gamma_4^{\frac{1}{N}}t^{\frac{1}{2}}\biggr[\log\biggr(e+\frac{1}{t}\biggr)\biggr]^\beta W_2(x,t)
 \int_0^t (t-s)^{-\frac{1}{2}}s^{-1}\biggr[\log\biggr(e+\frac{1}{s}\biggr)\biggr]^{-1-\beta}\,ds\\
 & \le C\gamma_4^{p-1}W_2(x,t)
\end{split}
\end{equation}
for $(x,t)\in D\times(0,1)$. 
Therefore, 
by \eqref{eq:5.12}, \eqref{eq:5.13} and \eqref{eq:5.20} we have 
\begin{equation*}
\begin{split}
 & \int_D G(x,y,t)\mu(y)\,dy+\int_0^t\int_{{\bf R}^{N-1}}G(x,y',0,t-s)\overline{u}(y',0,s)^p\,dy'\,ds\\
 & \le [S(t)\mu_1](x)+C\gamma_4^{p-1}W_1(x,t)
 +W_2(x,t)
 +C(cd)^p\gamma_4^{p-1} W_2(x,t)\\
 & \le[1+C\gamma_4^{p-1}][S(t)\mu_1](x)
 +c'[1+C(cd)^p\gamma_4^{p-1}]\Phi_\beta^{-1}\left([S(t)\Phi_\beta(\mu_2)](x)\right)
\end{split}
\end{equation*}
in $D\times(0,1)$ for some $c'>0$. 
Setting $d=2c'$ and taking a sufficiently small $\gamma_4>0$ if necessary, 
we obtain 
$$
\int_{{\bf R}^N_+} G(x,y,t)\mu(y)\,dy+\int_0^t\int_{{\bf R}^{N-1}}G(x,y',0,t-s)\overline{u}(y',0,s)^p\,dy'\,ds
\le\overline{u}(x,t)
$$
in $D\times(0,1)$. 
Therefore $\overline{u}$ is a supersolution of \eqref{eq:1.1} and \eqref{eq:1.2} in $[0,1)$, 
and by Lemma~\ref{Lemma:2.3} 
we can find a solution of \eqref{eq:1.1} and \eqref{eq:1.2} in $[0,1)$ such that 
$0\le u(x,t)\le\overline{u}(x,t)$ in $D\times(0,1)$. 
Thus Theorem~\ref{Theorem:1.5} follows. 
$\Box$
\section{Life span of solutions}
Since the minimal solution is unique, 
we can define the maximal existence time $T(\mu)$  of the minimal solution $u$ 
of \eqref{eq:1.1} and \eqref{eq:1.2}. 
We call $T(\mu)$ the life span of the solution~$u$. 
%
\subsection{Life span for large initial data}
Let $\kappa>0$ and $\varphi$ be a nonnegative measurable function in $D$.  
In this subsection we study the behavior of $T(\kappa\varphi)$ as $\kappa\to\infty$. 
\vspace{3pt}

Firstly, 
by Theorems~\ref{Theorem:1.1}, \ref{Theorem:1.3} and \ref{Theorem:1.4} 
we easily obtain the following result 
(compare with \cite[Theorem~5.1]{IS1} and \cite[Corollary~1.2]{IS2}). 
\begin{theorem}
\label{Theorem:6.1}
Let $p>1$ and $\varphi$ be a nonnegative continuous function in $D$ such that 
$$
0<\|\varphi\|_{L^\infty({\bf R}^{N-1}\times[0,\delta])}<\infty,
\qquad
\int_D e^{-\Lambda y_N^2}\varphi(y)\,dy<\infty,
$$
for some $\delta>0$ and $\Lambda>0$. 
Then there exists $\gamma>0$ such that 
$$
\gamma^{-1}\kappa^{-2(p-1)}\le T(\kappa\varphi)\le \gamma\kappa^{-2(p-1)}
$$
for all sufficiently large $\kappa>0$. 
\end{theorem}
Next we consider the case of $\mbox{dist}\,(\mbox{supp}\,\varphi,\partial D)>0$.
\begin{theorem}
\label{Theorem:6.2}
Let $p>1$ and $\varphi$ be a nonnegative measurable function in $D$ such that 
\begin{equation}
\label{eq:6.1}
\int_D e^{-\Lambda y_N^2}\varphi(y)\,dy<\infty
\end{equation}
for some $\Lambda\ge 0$. Assume that 
$L:=\mbox{dist}\,(\mbox{supp}\,\varphi,\partial D)>0$.  
Then 
\begin{equation}
\label{eq:6.2}
\lim_{\kappa\to\infty}\,(\log\kappa) T(\kappa\varphi)
=\frac{L^2}{4}. 
\end{equation}
\end{theorem}
{\bf Proof.}
We write $T_\kappa=T(\kappa\varphi)$ for simplicity. 
For any $\epsilon>0$, we can find $z=(z',z_N)\in D$ such that 
\begin{equation}
\label{eq:6.3}
\mbox{dist}\,(z,\partial D)\le L+\epsilon,
\qquad
\overline{\varphi}(z):=\lim_{r\to 0}\,\dashint_{B(z,r)}\varphi(y)\,dy>0. 
\end{equation}
Then, by Theorem~\ref{Theorem:1.1}, 
for any $\delta_1>0$, we can find $\gamma_1>0$ such that 
\begin{equation}
\label{eq:6.4}
\exp\left(-(1+\delta_1)\frac{z_N^2}{4T_\kappa}\right)
\,\dashint_{B_+(z,T_\kappa^\frac{1}{2})}\kappa\varphi(y)\,dy
\le \gamma_1\,T_\kappa^{-\frac{1}{2(p-1)}}.
\end{equation}
This implies that $T_\kappa\to 0$ as $\kappa\to\infty$. 
Furthermore, 
by \eqref{eq:6.3} and \eqref{eq:6.4} we have 
$$
\frac{\kappa}{2}\overline{\varphi}(z)\le\gamma_1 T_\kappa^{-\frac{1}{2(p-1)}}
\exp\left((1+\delta_1)\frac{z_N^2}{4T_\kappa}\right)
\le\exp\left((1+2\delta_1)\frac{(L+\epsilon)^2}{4T_\kappa}\right)
$$
for all sufficiently large $\kappa$. 
Then we obtain  
$$
(1-\delta_1)\log\kappa\le(1+2\delta_1)\frac{(L+\epsilon)^2}{4T_\kappa}
$$ 
for all sufficiently large $\kappa$, which implies that 
$$
\limsup_{\kappa\to\infty}\,(\log\kappa)T_\kappa\le\frac{1+2\delta_1}{1-\delta_1}\frac{(L+\epsilon)^2}{4}.
$$
Letting $\delta_1\to 0$ and $\epsilon\to 0$, we deduce that 
\begin{equation}
\label{eq:6.5}
\limsup_{\kappa\to\infty}\,
(\log\kappa)T_\kappa\le\frac{L^2}{4}.
\end{equation}

Let $0<d<L$.
Let $\delta_2\in(0,1)$ be such that 
\begin{equation}
\label{eq:6.6}
L>\frac{d}{(1-2\delta_2)^2}.
\end{equation} 
Set 
$\tilde{T}_\kappa:=d^2/4\log\kappa$ and $\lambda:=(1-\delta_2)/4\tilde{T}_\kappa$. 
Let $x\in D$ be such that $B(x,\tilde{T}_\kappa^{1/2})\,\cap\,\mbox{supp}\,\varphi\not=\emptyset$. 
Then, by \eqref{eq:6.1} we have
\begin{equation}
\label{eq:6.7}
\begin{split}
 & \int_{B(x,\tilde{T}_\kappa^{\frac{1}{2}})}
e^{-\lambda y_N^2}\varphi(y)\,dy
=\int_{B(x,\tilde{T}_\kappa^{\frac{1}{2}})}
e^{-(\lambda-\Lambda)y_N^2}e^{-\Lambda y_N^2}\varphi(y)\,dy\\
 & \le C\sup_{y\in B(x,\tilde{T}_\kappa^{\frac{1}{2}})}
\exp\left(-\frac{1-2\delta_2}{4\tilde{T}_\kappa}y_N^2\right)
\le C\exp\left(-\frac{(1-2\delta_2)\{L-2\tilde{T}_\kappa^{\frac{1}{2}}\}^2}{d^2}\log\kappa\right)\\
 & \le C\exp\left(-\frac{(1-2\delta_2)^2L^2}{d^2}\log\kappa\right)
 \le C\kappa^{-\frac{(1-2\delta_2)^2L^2}{d^2}}
\end{split}
\end{equation}
for all sufficiently large $\kappa$. 
This together with \eqref{eq:6.6} implies that 
$$
\sup_{x\in D}\,\dashint_{B(x,T^{\frac{1}{2}})}e^{-\lambda y_N^2}\kappa\varphi(y)\,dy
\le C\kappa^{1-\frac{1}{(1-2\delta_2)^2}}\tilde{T}_\kappa^{-\frac{N}{2}}
=o\left(\tilde{T}_\kappa^{-\frac{1}{p-1}}\right)
$$
as $\kappa\to\infty$. 
Therefore, by Theorem~\ref{Theorem:1.4} we see that $T_\kappa\ge\tilde{T}_\kappa$ 
for sufficiently large $\kappa$, and we obtain
$$
\liminf_{\kappa\to\infty}\,(\log\kappa)T_\kappa\ge\frac{d^2}{4}.
$$
Letting $d\to L$, we obtain 
$$
\liminf_{\kappa\to\infty}\,(\log\kappa)T_\kappa\ge\frac{L^2}{4}.
$$
This together with \eqref{eq:6.5} implies \eqref{eq:6.2}. 
Thus Theorem~\ref{Theorem:6.2} follows.
$\Box$
\vspace{5pt}

\noindent
Similarly, we have:
\begin{theorem}
\label{Theorem:6.3}
Let $p>1$, $z=(z',z_N)\in{\bf R}^N_+$ and $\delta_{z}(y):=\delta(y-z)$.   
Then 
$$
\lim_{\kappa\to\infty}\,(\log\kappa) T(\kappa\delta_z)
=\frac{z_N^2}{4}. 
$$
\end{theorem}

In the following two theorems, 
we study the relationship between 
the behavior of the life span $T(\kappa\varphi)$ 
for sufficiently large $\kappa$ and the singularity of $\varphi$ at $0\in\partial D$. 
Compare with \cite[Theorem~5.2]{IS1}. 
\begin{theorem}
\label{Theorem:6.4}
Let $\varphi$ be a nonnegative measurable function in $D$ such that 
\begin{equation}
\label{eq:6.8}
\varphi(y)\ge \gamma |y|^A\biggr[\log\left(e+\frac{1}{|y|}\right)\biggr]^{-B},
\qquad y\in B_+(0,1),
\end{equation}
for some $\gamma>0$, where $A>-N$, $B\in{\bf R}$ or $A=-N$, $B>1$. 
\vspace{3pt}
\newline
{\rm (i)} 
Let $1<p<p_*$. Then 
\begin{eqnarray}
\label{eq:6.9}
 & & 
\limsup_{\kappa\to\infty}
\frac{T(\kappa\varphi)}{[\kappa(\log\kappa)^{-B}]^{-\frac{2(p-1)}{A(p-1)+1}}}<\infty
\quad\qquad\mbox{if}\quad A>-N,\,\, B\in{\bf R},\vspace{7pt}\\
\label{eq:6.10}
 & &
\limsup_{\kappa\to\infty}
\frac{T(\kappa\varphi)}{[\kappa(\log\kappa)^{-B+1}]^{-\frac{2(p-1)}{A(p-1)+1}}}<\infty
\qquad\mbox{if}\quad A=-N,\,\, B>1. 
\end{eqnarray}
\newline
{\rm (ii)} Let $p>p_*$. If, either $A<-1/(p-1)$  and $B\in{\bf R}$ or $A=-1/(p-1)$ and $B<0$, 
then $T(\kappa\varphi)=0$ for all $\kappa>0$.  
If $A=-1/(p-1)$ and $B=0$, then $T(\kappa\varphi)=0$ for all sufficiently large $\kappa>0$. 
Furthermore, 
\begin{itemize}
  \item if $A>-1/(p-1)$, then \eqref{eq:6.9} holds:
  \item if $A=-1/(p-1)$ and $B>0$, then 
  $\displaystyle{\liminf_{\kappa\to\infty}\kappa^{-\frac{1}{B}}|\log T(\kappa\varphi)|>0}$. 
\end{itemize}
{\rm (iii)} Let $p=p_*$. If $A=-N$ and $B<N+1$, then 
$T(\kappa\varphi)=0$ for all $\kappa>0$. 
If $A=-N$ and $B=N+1$, then $T(\kappa\varphi)=0$ for all sufficiently large $\kappa>0$. 
Furthermore, 
\begin{itemize}
  \item if $A>-N$, then \eqref{eq:6.9} holds:
  \item if $A=-N$ and $B>N+1$, then 
  $\displaystyle{\liminf_{\kappa\to\infty}\kappa^{\frac{1}{-B+N+1}}|\log T(\kappa\varphi)|>0}$. 
\end{itemize}
\end{theorem}
{\bf Proof.}
We write $T_\kappa:=T(\kappa\varphi)$ for simplicity. 
We prove assertion~(i). 
Let $1<p<p_*$, $A>-N$ and $B\in{\bf R}$. 
For any $\epsilon\in(0,A+N)$, 
since 
$$
r^{-\epsilon}\biggr[\log\left(e+\frac{1}{r}\right)\biggr]^{-B}
\quad\mbox{is monotone decreasing near $r=0$},
$$
we have 
\begin{equation*}
\begin{split}
\int_{B(0,\sigma)}|y|^{A}\biggr[\log\left(e+\frac{1}{|y|}\right)\biggr]^{-B}\,dy
 & \ge C\sigma^{-\epsilon}\biggr[\log\left(e+\frac{1}{\sigma}\right)\biggr]^{-B}
\int_0^\sigma r^{A+\epsilon+N-1}\,dr\\
 & \ge C\sigma^{A+N}\biggr[\log\left(e+\frac{1}{\sigma}\right)\biggr]^{-B}
\end{split}
\end{equation*}
for all sufficiently small $\sigma>0$. 
Then 
it follows from Theorem~\ref{Theorem:1.1} that 
$T_\kappa\to 0$ as $\kappa\to\infty$ and 
$$
C\gamma\kappa T_\kappa^{\frac{A+N}{2}}\biggr[\log\left(e+T_\kappa^{-\frac{1}{2}}\right)\biggr]^{-B}
\le\gamma_1T_\kappa^{\frac{N}{2}-\frac{1}{2(p-1)}},
$$
that is
\begin{equation}
\label{eq:6.11}
T_\kappa^{\frac{A(p-1)+1}{2(p-1)}}\biggr[\log\left(e+T_\kappa^{-\frac{1}{2}}\right)\biggr]^{-B}
\le C\kappa^{-1}
\end{equation}
for all sufficiently large $\kappa>0$. 
Let $\tilde{T}_\kappa>0$ be such that 
$$
\tilde{T}_\kappa^{\frac{A}{2}+\frac{1}{2(p-1)}}\biggr[\log\left(e+\tilde{T}_\kappa^{-\frac{1}{2}}\right)\biggr]^{-B}
=C\kappa^{-1}. 
$$
Since $A>-N$ and $1<p<p_*$, 
$$
h(s):=s^{\frac{A(p-1)+1}{2(p-1)}}\biggr[\log\left(e+s^{-\frac{1}{2}}\right)\biggr]^{-B}
$$
is monotone increasing for all sufficiently small $s>0$.  
Then, by \eqref{eq:6.11} we have 
\begin{equation}
\label{eq:6.12}
T_\kappa\le\tilde{T}_\kappa
\le C[\kappa^{-1}(\log\kappa)^B]^{\frac{2(p-1)}{A(p-1)+1}}
\end{equation}
for all sufficiently large $\kappa>0$.  
This implies \eqref{eq:6.9}.

In the case where $A=-N$ and $B>1$, then 
\begin{equation*}
\begin{split}
\int_{B(0,\sigma)}|y|^{A}\biggr[\log\left(e+\frac{1}{|y|}\right)\biggr]^{-B}\,dy
 & \ge C\int_0^\sigma r^{-1}\biggr[\log\left(e+\frac{1}{r}\right)\biggr]^{-B}\,dr\\
 & \ge C\biggr[\log\left(e+\frac{1}{\sigma}\right)\biggr]^{-B+1}
\end{split}
\end{equation*}
for all sufficiently small $\sigma>0$. 
Then, similarly to \eqref{eq:6.11}, we obtain 
$$
T_\kappa^{\frac{A(p-1)+1}{2(p-1)}}\biggr[\log\left(e+T_\kappa^{-\frac{1}{2}}\right)\biggr]^{-B+1}
\le C\kappa^{-1}
$$
for all sufficiently large $\kappa>0$. 
Therefore, by a similar argument as in \eqref{eq:6.12} 
we obtain \eqref{eq:6.10}. Thus assertion~(i) follows. 

By similar arguments as in assertion~(i) 
we apply Theorem~\ref{Theorem:1.1} to obtain assertions~(ii) and (iii). 
(We leave the details of the proof to the reader.) 
Then Theorem~\ref{Theorem:6.4} follows. 
$\Box$
\begin{theorem}
\label{Theorem:6.5}
Let $\varphi$ be a nonnegative measurable function in $D$ such that 
$\mbox{supp}\,\varphi\subset B(0,1)$ and 
\begin{equation}
\label{eq:6.13}
\varphi(y)\le \gamma |y|^{A}\biggr[\log\left(e+\frac{1}{|y|}\right)\biggr]^{-B},
\qquad y\in B(0,1),
\end{equation}
for some $\gamma>0$. 
\vspace{3pt}
\newline
{\rm (i)} 
Let $1<p<p_*$. Then
\begin{eqnarray}
\label{eq:6.14}
 & & 
\liminf_{\kappa\to\infty}
\frac{T(\kappa\varphi)}{[\kappa(\log\kappa)^{-B}]^{-\frac{2(p-1)}{A(p-1)+1}}}>0
\quad\qquad\mbox{if}\quad A>-N,\,\, B\in{\bf R},\vspace{7pt}\\
\nonumber 
 & & 
\liminf_{\kappa\to\infty}
\frac{T(\kappa\varphi)}{[\kappa(\log\kappa)^{-B+1}]^{-\frac{2(p-1)}{A(p-1)+1}}}>0
\qquad\mbox{if}\quad A=-N,\,\, B>1. 
\end{eqnarray}
{\rm (ii)} Let $p>p_*$. 
\begin{itemize}
  \item If $A>-1/(p-1)$, then \eqref{eq:6.14} holds:
  \item If $A=-1/(p-1)$ and $B>0$, then 
  $\displaystyle{\limsup_{\kappa\to\infty}\kappa^{-\frac{1}{B}}|\log T(\kappa\varphi)| < \infty}$.
\end{itemize}
{\rm (iii)} Let $p=p_*$. 
\begin{itemize}
  \item If $A>-N$, then \eqref{eq:6.14} holds:
  \item If $A=-N$ and $B>N+1$, then 
  $\displaystyle{\limsup_{\kappa\to\infty}\kappa^{\frac{1}{-B+N+1}}|\log T(\kappa\varphi)|<\infty}$.
\end{itemize}
\end{theorem}
{\bf Proof.}
In the case where $A>-N$ and $B\in{\bf R}$,
for any $\epsilon\in(0,A+N)$, 
$$
r^{\epsilon}\biggr[\log\left(e+\frac{1}{r}\right)\biggr]^{-B}
\quad\mbox{is monotone increasing near $r=0$}
$$
and we have 
\begin{equation*}
\begin{split}
\int_{B(0,\sigma)}|y|^{A}\biggr[\log\left(e+\frac{1}{|y|}\right)\biggr]^{-B}\,dy
 & \le C\sigma^{\epsilon}\biggr[\log\left(e+\frac{1}{\sigma}\right)\biggr]^{-B}
\int_0^\sigma r^{A-\epsilon+N-1}\,dr\\
 & \le C\sigma^{A+N}\biggr[\log\left(e+\frac{1}{\sigma}\right)\biggr]^{-B}
\end{split}
\end{equation*}
for all sufficiently small $\sigma>0$. 
In the case where $A=-N$ and $B>1$, we obtain 
\begin{equation*}
\begin{split}
\int_{B(0,\sigma)}|y|^{A}\biggr[\log\left(e+\frac{1}{|y|}\right)\biggr]^{-B}\,dy
 & \le C\int_0^\sigma r^{-1}\biggr[\log\left(e+\frac{1}{r}\right)\biggr]^{-B}\,dr\\
 & \le C\biggr[\log\left(e+\frac{1}{\sigma}\right)\biggr]^{-B+1}
\end{split}
\end{equation*}
for all sufficiently small $\sigma>0$. 
Then Theorem~\ref{Theorem:6.5} follows from Theorems~\ref{Theorem:1.3}, 
\ref{Theorem:1.4} and  Theorem~\ref{Theorem:1.5}. 
We leave the details of the proof to the reader. 
(See also the proof of Theorem~\ref{Theorem:6.4}.)
$\Box$
\subsection{Life span for small initial data}
Motivated by \cite{LN}, 
we state two theorems on the behavior of $T(\kappa\varphi)$ as $\kappa\to 0$. 
Theorems~\ref{Theorem:6.6} and \ref{Theorem:6.7} follow 
from Theorem~\ref{Theorem:1.1} and Theorems~\ref{Theorem:1.3}--\ref{Theorem:1.5}. 
\begin{theorem}
\label{Theorem:6.6}
Let $N\ge1$ and $p>1$. 
Let $A>0$ and $\varphi$ be a nonnegative measurable function in $D$ such that 
$0\le\varphi(x)\le(1+|x|)^{-A}$ for $x\in D$.
\vspace{3pt}
\newline
{\rm (i)} 
Let $p=p_*$ and $A\ge 1/(p-1)=N$. 
Then there exists $\gamma>0$ such that 
$$
\log T(\kappa\varphi)\ge
\left\{
\begin{array}{ll}
  \gamma\kappa^{-(p-1)} & \mbox{if}\quad A>N,\\
  \gamma\kappa^{-\frac{p-1}{p}} & \mbox{if}\quad A=N,\\
\end{array}
\right.
$$
for all sufficiently small $\kappa>0$. 
\vspace{3pt}
\newline
{\rm (ii)}
Let $1<p<p_*$ or  $A<1/(p-1)$. 
Then there exists $\gamma'>0$ such that 
$$
T(\kappa\varphi)\ge
\left\{
\begin{array}{ll}
  \gamma'\kappa^{-\left(\frac{1}{2(p-1)}-\frac{1}{2} \min \{A,N\} \right)^{-1}} & \mbox{if}\quad A\not=N,\vspace{3pt}\\
  \gamma'\left(\frac{\kappa^{-1}}{\log(\kappa^{-1})}\right)^{\left(\frac{1}{2(p-1)}-\frac{N}{2}\right)^{-1}} & \mbox{if}\quad A=N,\\ \end{array}
\right.
$$
for all sufficiently small $\kappa>0$. 
\end{theorem}
\begin{theorem}
\label{Theorem:6.7}
Let $N\ge 1$ and $p>1$. 
Let $A>0$ and $\varphi$ be a nonnegative $L^\infty(D)$-function such that 
$\varphi(x)\ge(1+|x|)^{-A}$ for $x\in D$.
\vspace{3pt}
\newline
{\rm (i)}
Let $p=p_*$ and $A\ge 1/(p-1)=N$. 
Then there exists $\gamma>0$ such that 
\begin{equation}
\label{eq:6.15}
\log T(\kappa\varphi)\le
\left\{
\begin{array}{ll}
  \gamma\kappa^{-(p-1)} & \mbox{if}\quad A>N,\\
  \gamma\kappa^{-\frac{p-1}{p}} & \mbox{if}\quad A=N,\\
\end{array}
\right.
\end{equation}
for all sufficiently small $\kappa>0$. 
\vspace{3pt}
\newline
{\rm (ii)} 
Let $1<p<p_*$ or  $A<1/(p-1)$. 
Then there exists $\gamma'>0$ such that 
\begin{equation}
\label{eq:6.16}
T(\kappa\varphi)\le
\left\{
\begin{array}{ll}
  \gamma'\kappa^{-\left(\frac{1}{2(p-1)}-\frac{1}{2} \min \{A,N\} \right)^{-1}} & \mbox{if}\quad A\not=N,\vspace{3pt}\\
  \gamma'\left(\frac{\kappa^{-1}}{\log(\kappa^{-1})}\right)^{\left(\frac{1}{2(p-1)}-\frac{N}{2}\right)^{-1}} & \mbox{if}\quad A=N,\\
\end{array}
\right.
\end{equation}
for all sufficiently small $\kappa>0$. 
\end{theorem}
Theorems~\ref{Theorem:6.6} and \ref{Theorem:6.7} are proved by 
similar arguments as in \cite[Section~5]{HI}, 
and we leave the details of the proofs to the reader. 

Finally, we show 
that $\lim_{\kappa\to 0}T(\kappa\varphi)=\infty$ does not necessarily hold 
for problem~\eqref{eq:1.1}. 
\begin{theorem}
\label{Theorem:6.8}
Let $p>1$, $\lambda>0$ and $\varphi(x):=\exp{(\lambda x_N^2)}$.
Then 
\begin{equation}
\label{eq:6.17}
\lim_{\kappa\to+0} T(\kappa\varphi)=(4\lambda)^{-1}.
\end{equation}
\end{theorem}
{\bf Proof.}  Let $\kappa>0$ and $\delta>0$. 
Set $T_\kappa:= T(\kappa\varphi)$. 
It follows from Theorem~\ref{Theorem:1.1} that
\begin{equation*}
\begin{split}
\infty &> \exp\left(-(1+\delta)\frac{x_N^2}{4T_\kappa}\right)\int_{B(x,T_\kappa^\frac{1}{2})} \kappa\varphi(y)\,dy\\
& \ge C \exp\left(-(1+\delta)\frac{x_N^2}{4T_\kappa}\right) \kappa T_\kappa^{\frac{N}{2}}\exp\biggl(\lambda(x_N-T_\kappa^\frac{1}{2})^2\biggr)\\
& \ge C\kappa T_\kappa^{\frac{N}{2}} \exp\biggl\{ \left( \lambda-\frac{1+\delta}{4T_\kappa}\right)x_N^2\biggr\} 
\exp\biggl( -2\lambda T_\kappa^\frac{1}{2}x_N+\lambda T_\kappa \biggr)
\end{split}
\end{equation*}
for all $x\in D_{T_\kappa}$. 
This implies that $\lambda-(1+\delta)/4T_\kappa\le 0$. 
Since $\delta>0$ is arbitrary, we obtain 
\begin{equation} 
\label{eq:6.18}
\limsup_{\kappa\to +0}\,T_\kappa\le(4\lambda)^{-1}. 
\end{equation}
On the other hand, it follows that 
\begin{equation*}
\begin{split}
\dashint_{B(x, \tilde{T}_\delta^\frac{1}{2})}
\exp\left(-(1-\delta)\frac{y_N^2}{4\tilde{T}_\delta}\right)\kappa\exp{(\lambda y_N^2)}\,dy=\kappa,
\quad
x\in D_{\tilde{T}_\delta},
\end{split}
\end{equation*}
where $\tilde{T}_\delta := (1-\delta)/4\lambda$. 
Then we deduce from Theorem~\ref{Theorem:1.4} that $T_\kappa\ge \tilde{T}_\delta$ for all sufficiently small $\kappa>0$. 
Since $\delta>0$ is arbitrary, we obtain
$\liminf_{\kappa\to +0}T_\kappa\ge(4\lambda)^{-1}$. 
This together with \eqref{eq:6.18} implies \eqref{eq:6.17}. 
Thus Theorem~\ref{Theorem:6.8} follows.
$\Box$

\end{document}